\documentclass[12pt]{compositio}
\usepackage{amssymb}
\usepackage{amsmath}

\newcommand{\C}{\mathbb C}
\newcommand{\D}{\mathbb D}

\newcommand{\N}{\mathbb N}
\newcommand{\R}{\mathbb R}

\renewcommand{\P}{\mathbb P}

\newcommand{\Z}{\mathbb Z}
\newcommand{\cC}{\mathcal C}
\newcommand{\cI}{\mathcal I}
\newcommand{\cF}{\mathcal F}
\newcommand{\cM}{\mathcal M}

\newcommand{\re}{\mathrm{Re}}

\newcommand{\pa}{\partial}

\newcommand{\na}{\nabla}
\renewcommand{\a}{\alpha}
\renewcommand{\b}{\beta}
\newcommand{\g}{\gamma}

\newcommand{\De}{\mathit{\Delta}}
\renewcommand{\l}{\lambda}

\newcommand{\e}{\varepsilon}
\newcommand{\f}{\varphi}

\newcommand{\G}{\mathit{\Gamma}}

\newcommand{\W}{\mathit{\Omega}}
\newcommand{\w}{\omega}

\newcommand{\tr}{\;^t}

\newcommand{\diag}{\mathrm{diag}}
\newcommand{\X}{\frak X}

\newcommand{\bu}{\bullet}
\newcommand{\cH}{\mathcal{H}}
\newcommand{\cE}{\mathcal{E}}
\newcommand{\cO}{\mathcal{O}}

\newcommand{\pr}{\mathrm{pr}}
\newcommand{\id}{\mathrm{id}}

\newcommand{\qqed}{\qed \medskip\noindent }
\newtheorem{theorem}{Theorem}[section]

\newtheorem{proposition}{Proposition}[section]
\newtheorem{lemma}{Lemma}[section]
\newtheorem{cor}{Corollary}[section]
\newtheorem{fact}{Facts}[section]
\newtheorem{remark}{Remark}[section]

\makeatletter
\@addtoreset{equation}{section}
\def\theequation{\thesection.\arabic{equation}}
\makeatother
\title[Period relations for $F_4$]
{The monodromy representation and twisted period relations 
for Appell's hypergeometric function $F_4$}
\author{Yoshiaki Goto}
\email{y-goto@math.sci.hokudai.ac.jp }
\address[Goto]{
Department of Mathematics,
Hokkaido University,
Sapporo 060-0810, Japan
}

\author{Keiji Matsumoto}
\email{matsu@math.sci.hokudai.ac.jp}
\address[Matsumoto]{
Department of Mathematics,
Hokkaido University,
Sapporo 060-0810, Japan
}

\date{\today}
\dedication{Dedicated to Professor Kyoichi Takano on his seventieth birthday}

\begin{document}

\begin{abstract}
We consider the system $\cF_4(a,b,c)$ of differential equations 
annihilating Appell's hypergeometric series $F_4(a,b,c;x)$.
We find the integral representations 
for four linearly independent solutions expressed 
by the hypergeometric series $F_4$.
By using the intersection forms of twisted (co)homology groups 
associated with them, 
we provide the monodromy representation of $\cF_4(a,b,c)$ and the 
twisted period relations for the fundamental systems of 
solutions of $\cF_4$.
\end{abstract}
\maketitle

\section{Introduction}
Appell's hypergeometric series $F_4(a,b,c;x)$ 
of variables $x=(x_1,x_2)$ with complex parameters $a,b,c=(c_1,c_2)$ 
is defined by 
$$
F_4(a,b,c;x)=\sum_{(n_1,n_2)\in \N^2}
\frac{(a,n_1+n_2)(b,n_1+n_2)}
{(c_1,n_1)(c_2,n_2)(1,n_1)(1,n_2)}x_1^{n_1}x_2^{n_2},
$$
where $c_1,c_2\notin -\N=\{0,-1,-2,\dots\}$ 
and $(c_1,n_1)=c_1(c_1+1)\cdots(c_1+n-1)=\G(c_1+n_1)/\G(c_1)$.
This series converges in the set 
$$\D=\{x\in \C^2\mid \sqrt{|x_1|}+\sqrt{|x_2|}<1\},$$
satisfies 
$$F_4(a,b,c;x)=F_4(b,a,c;x),$$
and admits the integral representations (\ref{eq:K-rep}), 
(\ref{eq:A-rep}), and (\ref{eq:C-rep}).   
The system $\cF_4(a,b,c)$ of differential equations 
annihilating Appell's hypergeometric series $F_4(a,b,c;x)$ 
is a holonomic system of rank $4$ 
with the singular locus $S$ given in (\ref{eq:sing-loc}).
A fundamental system of solutions of $\cF_4(a,b,c)$ in 
a simply connected domain $U$ in $\D-S$
is expressed in terms of Appell's hypergeometric series $F_4$ with 
different parameters; see (\ref{eq:sols}) for 
their explicit forms.

In this paper, we find the twisted cycles associated with 
the integrand in (\ref{eq:K-rep}) which correspond to 
the solutions (\ref{eq:sols}). 
We evaluate the intersection numbers of several twisted cycles. 
By using the intersection numbers, as in \cite{M2} and \cite{MY},
we provide the monodromy representation  of $\cF_4(a,b,c)$; 
see Theorem \ref{th:monodromy}. 
We provide a basis for the twisted cohomology group associated with 
the integrand in (\ref{eq:K-rep}), and evaluate 
the intersection matrix for this basis; see Theorem \ref{th:int-c}.
By the compatibility of the parings of twisted (co)homology groups, 
we have the identity (\ref{eq:TPR}) 
for the intersection matrices 
and the period matrices for our bases of twisted (co)homology groups;  
for details, refer to Theorem \ref{th:TPR}.
This identity implies twisted period relations, which are quadratic relations 
between a fundamental system of solutions of $\cF_4$ 
and those of $\cF_4$ with different parameters. 
We present some examples in Corollary \ref{cor:TPR}.

There have been several studies of monodromy representations of the system $\cF_4(a,b,c)$
under the condition 
$$c_1,\ c_2,\ a,\ a-c_1,\ a-c_2,\ a-c_1-c_2,\ 
b,\ b-c_1,\ b-c_2,\ b-c_1-c_2\notin \Z;$$ 
see \cite{HU}, \cite{K}, and \cite{T}. 
It is determined in \cite{Ka2} that representation matrices are 
valid even when $c_1,c_2$ are positive integers,   
and that the system $\cF_4(a,b,c)$ is irreducible if and only if 
$c_1,c_2\notin \Z$ are removed from the above. 
Our expression of the monodromy representation is 
independent of the choice of fundamental systems of solutions of 
$\cF_4(a,b,c)$, and it is valid even in the case $c_1,c_2\in \Z$.
We represent circuit transforms as matrices by assigning fundamental systems 
of solutions of $\cF_4(a,b,c)$; 
see Corollary \ref{cor:mat-rep} and Remark  \ref{rem:mat-rep}.

Twisted period relations for Lauricella's system $\cF_D$ and 
Appell's system $\cF_2$, $\cF_3$ are studied in \cite{CM} and \cite{M1}.
We can obtain an explicit form of that for $\cF_4$ by  
evaluating the intersection matrix 
for the basis of the twisted cohomology group. 
We show that the intersection matrix $H$ of twisted cycles corresponding to 
the fundamental system of solutions of $\cF_4(a,b,c)$ in $U$ 
is diagonal. This fact is a key to obtaining several simple formulas for  
$F_4(a,b,c;x)$ that arise from the identity (\ref{eq:TPR}).
There is another application of the intersection form of twisted cohomology 
groups; we have a Pfaffian system of $\cF_4(a,b,c)$ using it 
as in \cite{M3}.  
For this, we refer the reader to the forthcoming paper \cite{GKM}.

Appell's system $\cF_4(a,b,c)$ is generalized to Lauricella's system 
$\cF_C(a,b,c)$ of rank $2^m$ with $m$-variables. 
A fundamental system of solutions of $\cF_C(a,b,c)$ 
near the origin is expressed 
in terms of  Lauricella's hypergeometric series $F_C(a,b,c;x)$. 
Their integral representations have been given in \cite{G1}; here, 
$2^m$ twisted cycles corresponding to them are constructed and  
the intersection numbers of these twisted cycles are evaluated. 
These results together with some intersection numbers of 
twisted closed $m$-forms imply that there are twisted period relations for the 
fundamental systems of $\cF_C$.
Similar results for Lauricella's system $\cF_A(a,b,c)$ have been obtained 
in \cite{G2}.

\section{Appell's system $\cF_4(a,b,c)$}
In this section, we collect some facts about 
Appell's system $\cF_4(a,b,c)$ of hypergeometric differential equations 
annihilating $F_4(a,b,c;x)$. 

Let $\pa_i$ $(i=1,2)$ be 
the partial differential operator with respect to $x_i$.
The function $F_4(a,b,c;x)$ satisfies differential equations
\begin{eqnarray*}
& &
\Big[x_1(1- x_1)\pa_1^2- x_2^2\pa_2^2- 2x_1x_2\pa_1\pa_2
 + \{c_1- (a+ b+ 1)x_1\}\pa_1- 
(a+ b+ 1)x_2\pa_2- ab\Big]f(x)=0,
\\
& &
\Big[x_2(1- x_2)\pa_2^2- x_1^2\pa_1^2- 2x_1x_2\pa_1\pa_2
 + \{c_2- (a+ b+ 1)x_2\}\pa_2- 
(a+ b+ 1)x_1\pa_1- ab\Big]f(x)=0.
\end{eqnarray*}
The system generated by them is called 
Appell's hypergeometric system $\mathcal{F}_4(a,b,c)$ 
of differential equations.
Though the function $F_4(a,b,c;x)$ is not defined for the case 
$c_1,c_2\in -\N$, the system $\mathcal{F}_4(a,b,c)$ is defined 
in this case, and
it is a holonomic system of rank $4$ with the singular locus 
\begin{equation}
\label{eq:sing-loc}
S=\{(x_1,x_2)\in \C^2\mid x_1x_2R(x)=0\}\cup L_\infty, \quad 
R(x)=x_1^2+x_2^2-2x_1x_2-2x_1-2x_2+1,
\end{equation}
where $L_\infty$ is the line at infinity in the projective plane $\P^2$.
We set $X=\P^2-S$.
We denote by $\cF_4(a,b,c;U)$ 
the vector space of solutions of $\mathcal{F}_4(a,b,c)$ in a 
simply connected domain $U\subset X\cap \D$.

If $c_1,c_2\notin \Z$, then $\cF_4(a,b,c;U)$ 
is spanned  by 
\begin{eqnarray}
\label{eq:sols}
& &\hspace{18mm}F_4(a,b,c;x),\\ 
\nonumber
& &\hspace{9mm}x_1^{1-c_1}F_4(a+1-c_1,b+1-c_1,2-c_1,c_2;x),\\
\nonumber
& &\hspace{9mm}x_2^{1-c_2}F_4(a+1-c_2,b+1-c_2,c_1,2-c_2;x),\\
\nonumber
& &x_1^{1-c_1}x_2^{1-c_2}F_4(a+2-c_1-c_2,b+2-c_1-c_2,2-c_1,2-c_2;x).
\end{eqnarray}
Note that $x_1^{1-c_1}$  and $x_2^{1-c_2}$ are single-valued 
holomorphic functions in $U$.

\begin{figure}[htb]
  \centering
\includegraphics[height=6cm,trim=0mm -8mm 0mm 0mm]{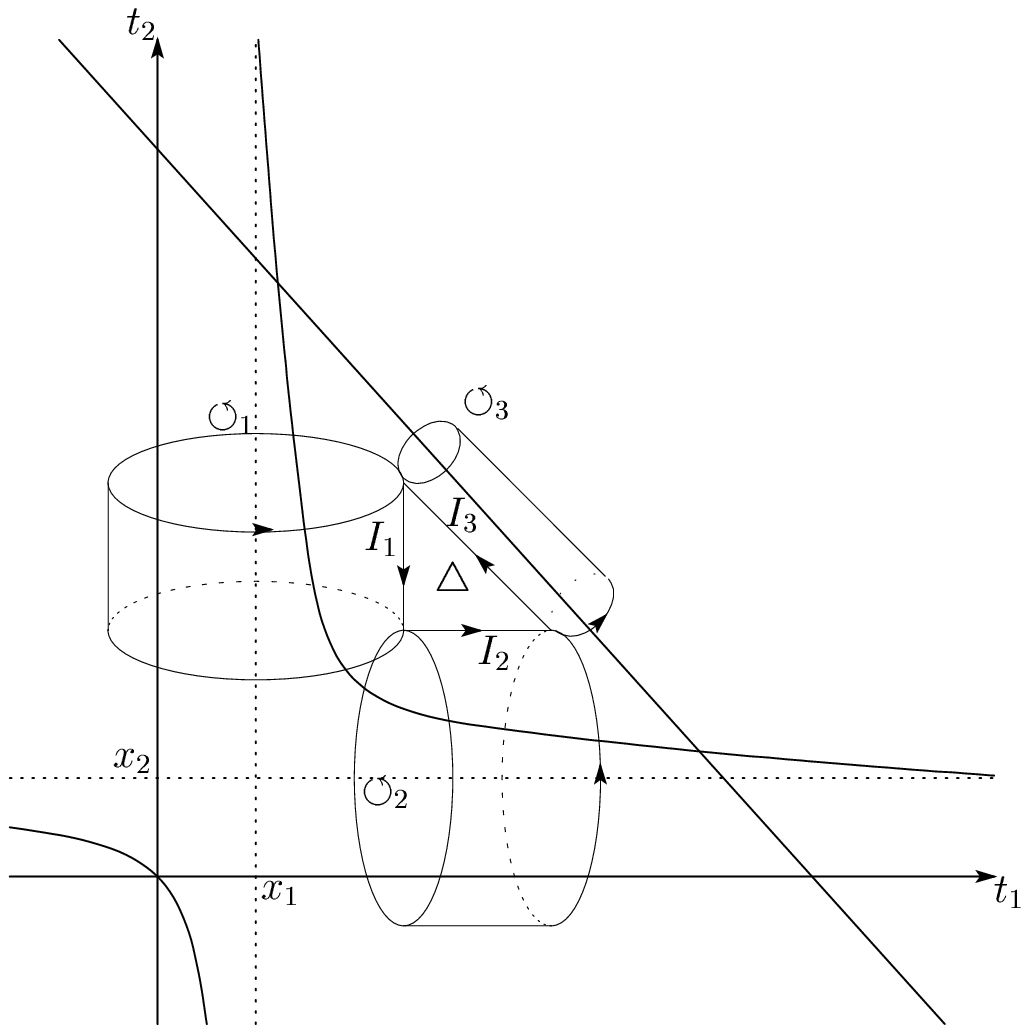}
\hspace{5mm}
\includegraphics[height=6cm,trim=0mm -5mm 0mm -3mm]{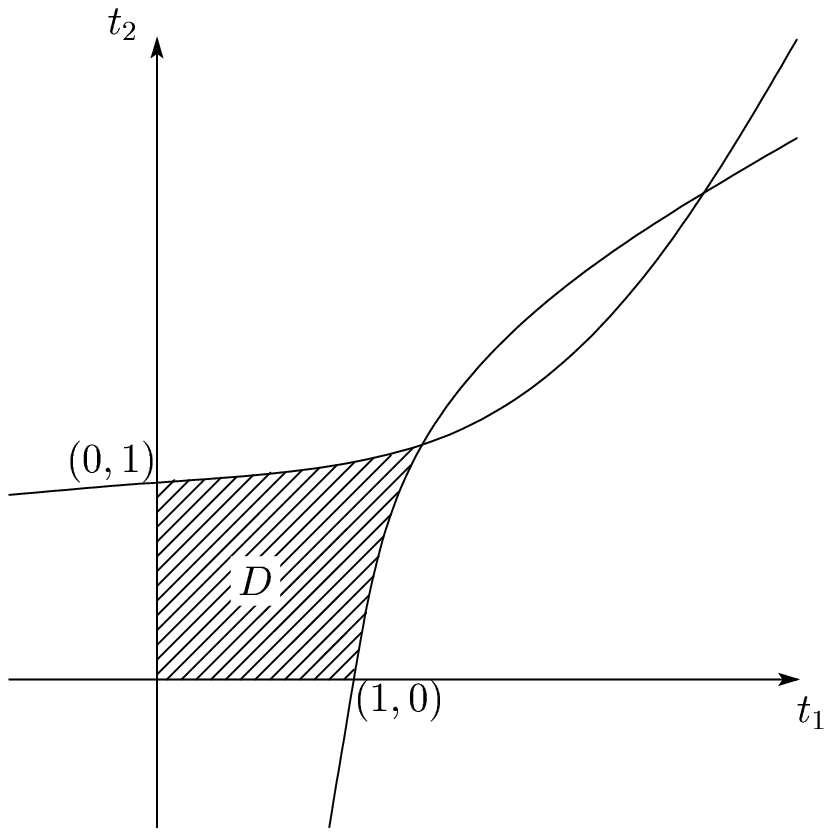}
  \caption{Domains of the integrals}
\label{fig:cycle}
\end{figure}

For sufficiently small positive real numbers $x_1$ and $x_2$, 
$F_4(a,b,c;x)$ admits
the following integral representations:
\begin{eqnarray}
\label{eq:K-rep}
& &
G_1
\int_{\De_1}t_1^{-c_1}t_2^{-c_2}(1\!-\! t_1\!-\! t_2)^{c_1+c_2-a-2}
\big(1-\frac{x_1}{t_1}-\frac{x_2}{t_2})^{-b}dt_1\wedge dt_2,\\
\nonumber
& &c_1,c_2,a-c_1-c_2\notin \Z,\\
\label{eq:A-rep}
& &G_2
\int_{\sqrt{-1}\R^2_x}
t_1^{-c_1}t_2^{-c_2}(1\!-\! t_1\!-\! t_2)^{c_1+c_2-a-2}
\big(1-\frac{x_1}{t_1}-\frac{x_2}{t_2})^{-b}dt_1\wedge dt_2,\\
\nonumber
& &\re(c_1-a)<1,\ \re(c_2-a)<1,\\
\label{eq:C-rep}
& &G_3
\int_{D}
t_1^{a-1}t_2^{b-1}(1\!-\! t_1\!+\! t_1t_2x_2)^{c_1-a-1}
(1\!-\! t_2\!+\! t_1t_2x_1)^{c_2-b-1}dt_1\wedge dt_2,
\\
\nonumber 
& &\re(c_1)>\re(a)>0,\ \re(c_2)>\re(b)>0.\ 
\end{eqnarray}
Here 
\begin{eqnarray*}
G_1&=&\frac{\G(1-a)}{\G(1-c_1)\G(1-c_2)\G(c_1+c_2-a-1)},\\
G_2&=&\frac{\G(c_1)\G(c_2)\G(a-c_1-c_2+2)}{(2\pi\sqrt{-1})^2\G(a)},\\
G_3&=&\frac{\G(c_1)\G(c_2)}{\G(a)\G(b)\G(c_1-a)\G(c_2-b)},
\end{eqnarray*}
$\De_1$ is the formal sum 
\begin{eqnarray*}
\De_1&=&\triangle+\frac{(\circlearrowleft_1\times I_1)}{1-\g_1^{-1}}
+\frac{(\circlearrowleft_2\times I_2)}{1-\g_2^{-1}}
+\frac{(\circlearrowleft_3\times I_3)}{1-\g_1\g_2\a^{-1}}\\
& &+\frac{(\circlearrowleft_1\times \circlearrowleft_2)}{(1-\g_1^{-1})(1-\g_2^{-1})}
+\frac{(\circlearrowleft_2\times \circlearrowleft_3)}{(1-\g_2^{-1})(1-\g_1\g_2\a^{-1})}
+\frac{(\circlearrowleft_3\times \circlearrowleft_1)}{(1-\g_1\g_2\a^{-1})(1-\g_1^{-1})},
\end{eqnarray*}
of $2$-dimensional real surfaces, $\triangle$ and its boundary components $I_i$
$(i=1,2,3)$ are given in Figure \ref{fig:cycle}, 
$\circlearrowleft_i$ $(i=1,2)$ is a positively oriented circle in the $t_i$-space 
starting from the projection of $I_i$ to this space and 
surrounding the divisors $t_i=0$, and $Q(t,x)=t_1t_2-t_1x_2-t_2x_1=0$
for $t\in I_i$, 
$\circlearrowleft_3$ is a positively oriented circle with a small radius 
in the orthogonal complement of the divisor $L(t)=1-t_1-t_2=0$ starting 
from the projection of $I_3$ to this space and surrounding the divisor,
$\a=e^{2\pi\sqrt{-1}a}$, 
$\b=e^{2\pi\sqrt{-1}b}$, $\g_i=e^{2\pi\sqrt{-1}c_i}$ $(i=1,2)$, 
$$\sqrt{-1}\R^2_x=\{(\sqrt{x_1},\sqrt{x_2})+
(s_1,s_2)\sqrt{-1}\mid s_1,s_2\in \R\}\subset \C^2,
\quad (\sqrt{x_1},\sqrt{x_2})\in \triangle,$$
and $D$ is the bounded connected component of 
$$\{(t_1,t_2)\in \R^2\mid t_1,t_2,1-t_1+t_1t_2x_2,1-t_2+t_1t_2x_1>0\};$$
see Figure \ref{fig:cycle}. 
The argument of each factor of the integrand of (\ref{eq:K-rep}) 
at any point $t=(t_1,t_2)\in \triangle$ 
is $0$, 
that of (\ref{eq:K-rep})  at the starting point of the circle 
$\circlearrowleft_i$ $(i=1,2,3)$ is $0$, 
that of (\ref{eq:A-rep}) at $(t_1,t_2)=(\sqrt{x_1},\sqrt{x_2})$  
is $0$, and that of (\ref{eq:C-rep}) 
at any point $t=(t_1,t_2)\in D$ is $0$.
For these integral representations 
of $F_4(a,b,c;x)$, we 
refer the reader to \cite{AoKi}, \cite{O}, and \cite{C}.

\begin{figure}[htb]
\centering
\begin{minipage}[h][8cm][t]{15cm}
\includegraphics[width=7cm]{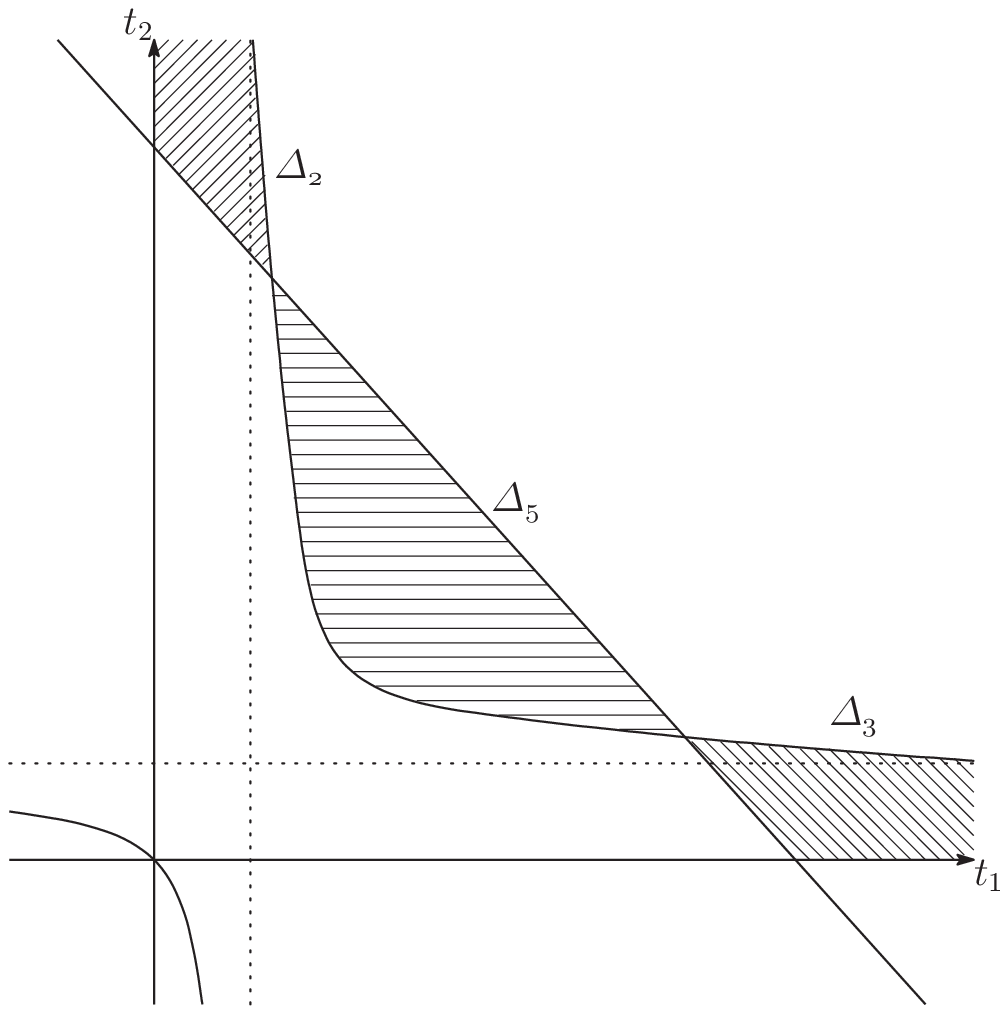}\hspace{5mm}
\raisebox{3.5cm}{
\begin{tabular}[htb]{c}
The arguments of the factors of the integrand\\[2mm]
$\displaystyle{\begin{array}{|c|c|c|c|c|}
\hline
&t_1& t_2 & 1-t_1-t_2& 1-\dfrac{x_1}{t_1}-\dfrac{x_2}{t_2}\\
\hline
\De_2& 0& 0& -\pi& -\pi\\
\hline
\De_3& 0& 0& -\pi& -\pi\\
\hline
\De_5& 0& 0& 0& 0\\
\hline
\end{array}}
$\\[2mm]
\end{tabular}
}
\end{minipage}
\label{fig:domain}
\caption{Domains of integrals}
\end{figure}

For $x\in U$, we set 
\begin{equation}
\label{eq:intrep}
f_i(x)=\int_{\De_i}t_1^{-c_1}t_2^{-c_2}(1-t_1-t_2)^{c_1+c_2-a-2} 
\big(1-\frac{x_1}{t_1}-\frac{x_2}{t_2}\big)^{-b}dt_1\wedge dt_2
,\quad (i=1,\dots,5),
\end{equation}
where 
$\De_2$, $\De_3$, and $\De_5$ are given in Figure \ref{fig:domain}, and 
$\De_4$ is the image of $\De_1$ under the involution 
$$\imath:(t_1,t_2)\mapsto (\frac{x_1}{t_1},\frac{x_2}{t_2}),$$
on $$\C_x^2=\{(t_1,t_2)\in \C^2
\mid t_1t_2(1-t_1-t_2)(t_1t_2-t_1x_2-t_2x_1)\ne0\}.
$$

\noindent
The conditions for their convergence are as follows.
\begin{table}[hbt]
$$
\begin{array}{|c|c|}
\hline
f_1&  c_1,c_2,a-c_1-c_2\notin \Z\\
\hline
f_2&  \re(b-c_1+1),\re(c_1+c_2-a-1),\re(1-b),\re(a-c_1+1)>0\\
\hline
f_3&  \re(b-c_2+1),\re(c_1+c_2-a-1),\re(1-b),\re(a-c_2+1)>0\\
\hline
f_4&  c_1,c_2,b-c_1-c_2\notin \Z\\
\hline
f_5&  \re(c_1+c_2-a-1),\re(1-b)>0\\
\hline
\end{array}
$$
\label{tab:converge}
\caption{Convergence conditions}
\end{table}

\begin{lemma}
\label{lem:series-rep}
We have
\begin{eqnarray*}
f_1(x)&=&\frac{\G(1-c_1)\G(1-c_2)\G(c_1+c_2-a-1)}{\G(1-a)}F_4(a,b,c_1,c_2;x),
\\[3mm]
f_2(x)&=&
\frac{\G(a+1-c_1)\G(b+1-c_1)\G(1-b)\G(c_1+c_2-a-1)}{\G(2-c_1)\G(c_2)}\\
& &\times e^{-\pi\sqrt{-1}(c_1+c_2-a-b)}x_1^{1-c_1}
F_4(a+1-c_1,b+1-c_1,2-c_1,c_2;x),\\[3mm]
f_3(x)&=&
\frac{\G(a+1-c_2)\G(b+1-c_2)\G(1-b)\G(c_1+c_2-a-1)}{\G(c_1)\G(2-c_2)}\\
& &\times e^{-\pi\sqrt{-1}(c_1+c_2-a-b)}
x_2^{1-c_2}F_4(a+1-c_2,b+1-c_2,c_1,2-c_2;x),\\[3mm]
f_4(x)&=&\frac{\G(c_1-1)\G(c_2-1)\G(1-b)}{\G(c_1+c_2-b-1)}
x_1^{1-c_1}x_2^{1-c_2}F_4(a\!+\!2\!-\! c_1\!-\! c_2,b\!+\!2\!-\! c_1\!-\! c_2,
2\!-\! c_1,2\!-\! c_2;x).
\end{eqnarray*}
\end{lemma}
\proof
Note that the first equality is nothing but 
the integral representation (\ref{eq:K-rep}).
We will show the last equality.
The transformation $\imath$ satisfies $\imath=\imath^{-1}$, and 
it implies  
\begin{eqnarray*}
f_4&=&x_1^{1-c_1}x_2^{1-c_2}\!\!
\int_{\De_1}t_1^{c_1-2}t_2^{c_2-2}
\big(1-\frac{x_1}{t_1}-\frac{x_2}{t_2}\big)^{c_1+ c_2- a-2}
(1- t_1- t_2)^{-b}dt_1\wedge dt_2\\
&=& x_1^{1-c_1}x_2^{1-c_2}
\frac{\G(c_1-1)\G(c_2-1)\G(1-b)}{\G(c_1+c_2-b-1)}
F_4(b\!+\! 2\!-\!c_1\!-\!c_2,a\!+\!2\!-\! c_1\!-\! c_2,
2\!-\!c_1,2\!-\!c_2;x).
\end{eqnarray*}
To obtain the second equality, we use an orientation-reversing transformation 
$$(s_1,s_2)\mapsto (t_1,t_2)=\big(x_1s_1,\frac{1}{s_2}\big),$$ 
which sends the domain $D$ to $\De_2$.
This transformation leads to 
\begin{eqnarray*}
f_2&=&-x_1^{1-c_1}\!\!
\int_{-D}s_1^{-c_1}s_2^{c_2-2}
\big(1\!-\! x_1s_1\!-\!\frac{1}{s_2}\big)^{c_1\!+\! c_2\!-\! a\!-\!2}
\big(1\!-\! \frac{1}{s_1}\!-\!s_2x_2\big)^{-b}ds_1\wedge ds_2\\
&=&x_1^{1-c_1}\!\!
\int_{D}s_1^{b-c_1}s_2^{a\!-\! c_1}
(s_2\!-\! x_1s_1s_2\!-\!1)^{c_1\!+\! c_2\!-\! a\!-\!2}
(s_1\!-\!1- x_2s_1s_2)^{-b}ds_1\wedge ds_2\\
&=&e^{-\pi\sqrt{-1}(c_1+c_2-a-b)}
x_1^{1-c_1}
\frac{\G(b+1-c_1)\G(a+1-c_1)\G(1-b)\G(c_1+c_2-a-1)}{\G(2-c_1)\G(c_2)}\\
& &\times F_4(b\!+\!1\!-\! c_1,a\!+\!1\!-\! c_1,2\!-\!c_1,c_2;x)
\end{eqnarray*}
by (\ref{eq:C-rep}). We can obtain the third equality in a similar way.
\qqed

\section{Twisted homology group}
Below, 
we will regard the parameters $a$, $b$, $c_1$, and $c_2$ as indeterminants, 
and we will assume that 
\begin{equation}
\label{eq:non-integral}
a,\  a-c_1,\ a-c_2,\  a-c_1-c_2,\  b,\  b-c_1,\ b-c_2,\  b-c_1-c_2,\ c_1,\ c_2
\notin \Z, 
\end{equation}
when we assign them to complex numbers.
Set 
$$\l_1=b-c_1+1,\quad \l_2=b-c_2+1,\quad 
\l_3=c_1+c_2-a-1,\quad 
\l_4=-b,
$$
and let $\C(\mu)$ be the rational function field of 
$\mu_1=e^{2\pi\sqrt{-1}\l_1},\dots,\mu_4=e^{2\pi\sqrt{-1}\l_4}$ 
over $\C$.

We define a subset $\X$ in $(\P^1\times \P^1)\times \P^2$ by 
$$\X=\{(t,x)\in \C^2\times X\mid
t_1t_2L(t)Q(t,x)\ne0\},\quad 
L(t)=1-t_1-t_2,\ Q(t,x)=t_1t_2-t_2x_1-t_1x_2.
$$
There is a natural projection 
$$\pr:\X\ni (t,x)\mapsto x\in X;$$
note that $\C_x^2=\pr^{-1}(x)$ for a fixed $x\in X$.
Let 
$$u=u(t,x)=t_1^{\l_1}t_2^{\l_2}L(t)^{\l_3}Q(t,x)^{\l_4}
=t_1^{b+1-c_1}t_2^{b+1-c_2}L(t)^{c_1+c_2-a-1}Q(t,x)^{-b}
$$
be a function of $(t,x)$ in a simply connected neighborhood 
of $(\dot t,\dot x)=\frac{1}{8}(\sqrt2,\sqrt2 ,1,1)\in \X$.
Along any path in $\X$ starting with $(\dot t,\dot x)$,  
we can make the analytic continuation of $u$. 
Though this continuation depends on the path, 
it is single valued and holomorphic around the end point of the path.

Let $\sigma$ be a $k$-chain in $\C_{x}^2$ for a fixed $x\in X$. 
We define a twisted $k$-chain $\sigma^u $ by 
$\sigma$ loading a branch of $u$ on it.
We denote the $\C(\mu)$-vector space of finite sums of twisted $k$-chains
by $\cC_k(\C_x^2,u)$. 
We define the boundary operator $\pa^u:\cC_k(\C_x^2,u)\to \cC_{k-1}(\C_x^2,u)$ 
by 
$$\sigma^u \mapsto \pa(\sigma)^{u|_{\pa(\sigma)}},$$
where $\pa$ is the usual boundary operator and 
$u|_{\pa(\sigma)}$ is the restriction of $u$ to $\pa(\sigma)$. 
We have a complex 
$$\cC_\bu(\C_x^2,u):\cdots \overset{\pa^u}{\longrightarrow} \cC_k(\C_x^2,u)
\overset{\pa^u}{\longrightarrow}  \cC_{k-1}(\C_x^2,u)
\overset{\pa^u}{\longrightarrow}\cdots,
$$
and its $k$-th homology group $H_k(\cC_\bu(\C_x^2,u))$.
Similarly we have 
a complex $\cC_\bu^{lf}(\C_x^2,u)$ of 
locally finite sums of twisted chains 
and its $k$-th homology group $H_k(\cC_\bu^{lf}(\C_x^2,u))$.
It is shown in \cite{AoKi} that 
$$H_k(\cC_\bu(\C_x^2,u))\simeq H_k(\cC_\bu^{lf}(\C_x^2,u)),
\quad \dim_{\C(\mu)}  H_k(\cC_\bu(\C_x^2,u))=
\left\{
\begin{array}{cl}
4 &\textrm{if }  k=2,\\
0 &\textrm{otherwise,} 
\end{array}
\right.
$$
for any fixed $x\in X$. Thus we have a map 
$$\mathrm{reg}:H_2(\cC_\bu^{lf}(\C_x^2,u))\to H_2(\cC_\bu(\C_x^2,u)),$$
which is the inverse of the natural map $ 
H_2(\cC_\bu(\C_x^2,u))\to H_2(\cC_\bu^{lf}(\C_x^2,u))$.

We regard the integral (\ref{eq:intrep}) as the pairing between 
the form 
$$\f_1=d\log\big(\frac{t_1}{L(t)}\big)\wedge d\log\big(\frac{t_2}{L(t)}\big)
=\frac{dt_1\wedge dt_2}{t_1t_2L(t)}$$
and $\De_i$ loaded with a branch of $u$, which represents 
an element of 
$H_2(\cC_\bu^{lf} (\C_x^2,u))$ $(i=1,\dots,5)$.  
The images of the element above under the map $\mathrm{reg}$ 
will be denoted by $\De_i^u\in H_2(\cC_\bu(\C_x^2,u))$ for $i=1,\dots,5$. 

By considering $1/u$ instead of $u$, we have 
$H_2(\cC_\bu(\C_x^2,1/u))$ and its elements 
$\De_1^{1/u}$,\dots,$\De_5^{1/u}$.
There is the intersection pairing $\cI_h$ between 
$H_2(\cC_\bu(\C_x^2,u))$ and $H_2(\cC_\bu(\C_x^2,1/u))$. 
It is defined as follows. 
Let $\De^u$ and $\Acute{\De}^{1/u}$ be elements of 
$H_2(\cC_\bu(\C_x^2,u))$ and $H_2(\cC_\bu(\C_x^2,1/u))$ 
given by 
$$
\De^u=\sum_{i\in I}  d_i D_i^{u_i},\quad 
\Acute{\De}^{1/u}=\sum_{j\in J}  \Acute{d}_j  \Acute{D}_j^{1/u_j}, \qquad d_i,
\Acute{d}_j\in\C(\mu),$$
where $D_i^{u_i}$ denotes a singular 
$2$-simplex $D_i$ loaded with  a branch $u_i$ of $u$.
Then their intersection number is 
$$
\cI_h(\De^u,\Acute{\De}^{1/u})
=\sum_{i\in I,j\in J}\sum_{p\in D_i\cap \Acute{D}_j}
  d_i \Acute{d}_j (D_i\cdot \Acute{D}_j)_p\frac{u_i(p)}{u_j(p)},$$
where 
$(D_i\cdot \Acute{D}_j)_p$ is the topological intersection number 
of $2$-chains $D_i$ and $\Acute{D_j}$ at $p$.
The intersection from $\cI_h$ is bilinear.
Since 
$$\De^{1/u}=\sum_{i\in I}  d_i^\vee D_i^{1/u_i},\quad 
\Acute{\De}^{u}=\sum_{j\in J}  \Acute{d}_j^\vee  \Acute{D}_j^{u_j},
$$
for the above $\De^u$ and $\Acute{\De}^{1/u}$, 
we have 
\begin{equation}
\label{eq:transpose}
\cI_h(\Acute{\De}^u,\De^{1/u})
=\cI_h(\De^u,\Acute{\De}^{1/u})^\vee,
\end{equation}
where $z(\mu_1,\dots,\mu_4)^\vee =z(1/\mu_1,\dots,1/\mu_4)$
for $z(\mu_1,\dots,\mu_4)\in \C(\mu)$.

\begin{lemma}
\label{lem:H-int}
The intersection numbers 
$\cI_h(\De_i^u ,\De_i^{1/u})$ ($i=1,\dots,4$) are 
\begin{eqnarray*}
\cI_h(\De_1^u ,\De_1^{1/u})
&=&\dfrac{1-(\mu_1\mu_4)(\mu_2\mu_4)(\mu_3)}
{(1-\mu_1\mu_4)(1-\mu_2\mu_4)(1-\mu_3)}=
\dfrac{-(1-\a)\g_1\g_2}{(\a-\g_1\g_2)(1-\g_1)(1-\g_2)}, \\
\cI_h(\De_2^u ,\De_2^{1/u})
&=&
\dfrac{(1-\mu_1\mu_4)(1-\mu_3(\mu_2\mu_3\mu_4)^{-1})}
{(1\!-\!\mu_1)(1\!-\!\mu_4)(1\!-\!\mu_3)(1\!-\!(\mu_2\mu_3\mu_4)^{-1})}=
\dfrac{\a\b\g_1(1-\g_1)(1-\g_2)}
{(\a\!-\!\g_1)(\a\!-\!\g_1\g_2)(\b\!-\!\g_1)(1\!-\!\b)}, \\
\cI_h(\De_3^u ,\De_3^{1/u})
&=&
\dfrac{(1-\mu_2\mu_4)(1-\mu_3(\mu_1\mu_3\mu_4)^{-1})}
{(1\!-\!\mu_2)(1\!-\!\mu_4)(1\!-\!\mu_3)(1\!-\!(\mu_1\mu_3\mu_4)^{-1})}=
\dfrac{\a\b\g_2(1-\g_1)(1-\g_2)}
{(\a\!-\!\g_2)(\a\!-\!\g_1\g_2)(1\!-\!\b)(\b\!-\!\g_2)}, \\
\cI_h(\De_4^u ,\De_4^{1/u})
&=&
\dfrac{1-(\mu_1\mu_4)^{-1}(\mu_2\mu_4)^{-1}(\mu_4)}
{(1-(\mu_1\mu_4)^{-1})(1-(\mu_2\mu_4)^{-1})(1-\mu_4)}=
\dfrac{-(\b-\g_1\g_2)}{(1-\b)(1-\g_1)(1-\g_2)}.
\end{eqnarray*}
\end{lemma}
\proof
To compute $\cI_h(\De_1^u ,\De_1^{1/u})$, we have only to 
follow \textsc{Example} 3.1 in Section 3 of Chapter VIII of \cite{Y2},
by considering the contribution of the divisor $Q(t,x)=0$. 
By using the involution $\imath$, we can evaluate 
$\cI_h(\De_4^u ,\De_4^{1/u})$.
For the rest, 
transform $\De_i$ $(i=2,3)$ to the domain $D$ in 
the expression (\ref{eq:C-rep}) as in the proof of Lemma \ref{lem:series-rep}; 
regard it as a quadrilateral and 
apply \textsc{Example} 3.2 in Section 3 of Chapter VIII of \cite{Y2}. 
\qqed

For a small simply connected neighborhood $U$ of $\dot x$, 
we have a family 
$$\bigcup_{x\in U}H_2(\cC_\bu(\C_x^2,u)),$$
which can be naturally identified with $\cF_4(a,b,c;U)$ by (\ref{eq:intrep}).
Since a path $\rho_x$ in $X$ connecting $\dot x$ and $x$ defines
the isomorphism 
$$(\rho_x)_*:H_2(\cC_\bu(\C_{\dot x}^2,u))\to H_2(\cC_\bu(\C_x^2,u)),$$   
we have a local system
$$\cH_2(X)=\bigcup_{x\in X}H_2(\cC_\bu(\C_x^2,u))$$
over $X$. 
Its stalk over $x$ is denoted by $H_2(\cC_\bu(\C_x^2,u))$.

Similarly, we have a local system 
$$\cH_2^\vee (X)=\bigcup_{x\in X}H_2(\cC_\bu(\C_x^2,1/u))$$
over $X$ with respect to $1/u$. The local triviality of these local systems
$\cH_2(X)$ and $\cH_2^\vee (X)$ 
imply the following.
\begin{proposition}
\label{prop:inv-h-int}
The intersection number is invariant under the deformation, that is, 
$$\cI_h((\rho_x)_*(\De^u),(\rho_x)_*(\Acute{\De}^{1/u}))
=\cI_h(\De^u,\Acute{\De}^{1/u})$$
for any $\De^u\in H_2(\cC_\bu(\C_{\dot x}^2,u))$, 
$\Acute{\De}^{1/u}\in H_2(\cC_\bu(\C_{\dot x}^2,1/u))$, 
and any path $\rho_x$ in $X$ connecting $\dot x$ and $x$.
\end{proposition}

\section{Monodromy representation}
A loop $\rho$ in $X$ with base point $\dot x$ 
induces a linear transformation of 
the stalk $H_2(\cC_\bu(\C_{\dot x}^2,u))$ of $\cH_2(X)$ over $\dot x$.
By this correspondence, we have a homomorphism 
$$\cM:\pi_1(X,\dot x)\to GL(H_2(\cC_\bu(\C_{\dot x}^2,u))),$$
which is called the monodromy representation of the local system $\cH_2(X)$.
Note that we can regard it as the monodromy representation of 
the system $\cF_4(a,b,c)$ by the 
identification of $\cF_4(a,b,c;U)$ for a small neighborhood $U$ of $\dot x$
with $\bigcup_{x\in U}H_2(\cC_\bu(\C_x^2,u))$.
It is shown in \cite{K} that 
the fundamental group $\pi_1(X,\dot x)$ is generated by three loops 
$\rho_i:[0,1] \to X$ ($i=1,2,3$), 
\begin{eqnarray*}
\rho_1&:&\theta \mapsto 
\left(\frac{\exp(2\pi\sqrt{-1}\theta)}{8},\frac{1}{8}\right),\\
\rho_2&:&\theta \mapsto 
\left(\frac{1}{8},\frac{\exp(2\pi\sqrt{-1}\theta)}{8}\right),\\
\rho_3&:&\theta \mapsto \left(\frac{2-\exp(2\pi\sqrt{-1}\theta)}{8},
\frac{2-\exp(2\pi\sqrt{-1}\theta)}{8}\right).
\end{eqnarray*}
Note that the loop $\rho_i$ $(i=1,2)$ turns the divisor $x_i=0$ 
positively,
and $\rho_3$ turns the divisor $R(x)=0$ positively.
We put $\cM_i=\cM(\rho_i)$ $(i=1,2,3)$.

\begin{proposition}
\label{prop:M1M2}
The elements $\De_1^u,\dots,\De_4^u$ span $H_2(\cC_\bu(\C_{\dot x}^2,u))$. 
With respect to the basis $\tr(\De_1^u,\dots,\De_4^u)$, 
$\cM_1$ and $\cM_2$ are represented by matrices 
$$\diag(1,\g_1^{-1},1,\g_1^{-1})\quad\textrm{and}\quad 
\diag(1,1,\g_2^{-1},\g_2^{-1}),$$
respectively, where $\diag(z_1,\dots,z_n)$ denotes the diagonal matrix with 
diagonal entries $z_1,\dots,z_n$.
\end{proposition}
\proof
Recall that the solutions $f_i$ are defined by the integrals over $\De_i$ 
in (\ref{eq:intrep}), and that they admit local expressions as in Lemma 
\ref{lem:series-rep}. 
We have
$$\tr\left((\rho_1)_*(\De_1^u),\dots,(\rho_1)_*(\De_4^u)\right)
=\diag(1,\g_1^{-1},1,\g_1^{-1})
\tr\left(\De_1^u,\dots,\De_4^u\right),$$
$$\tr\left((\rho_2)_*(\De_1^u),\dots,(\rho_2)_*(\De_4^u)\right)
=\diag(1,1,\g_2^{-1},\g_2^{-1})
\tr\left(\De_1^u,\dots,\De_4^u\right),$$
since the local behavior of $f_i$ is same to that of $\De_i$. 
\qqed

\begin{lemma}
\label{lem:intH}
If $i\ne j$ ($1\le i,j\le 4$) then 
$$\cI_h(\De_i^u,\De_j^{1/u})=0.$$ 
The intersection matrix 
$H=\left(\cI_h(\De_i^u,\De_j^{1/u})\right)_{1\le i,j\le4}$
is a diagonal  matrix with entries as given in Lemma \ref{lem:H-int}.
\end{lemma}
\proof
By Propositions \ref{prop:inv-h-int} and \ref{prop:M1M2}, we have 
$$\cI_h(\De_i^u,\De_j^{1/u})=
\cI_h((\rho_1)_*(\De_i^u),(\rho_1)_*(\De_j^{1/u}))
=\cI_h(\g_1^{-1}\De_i^u,\De_j^{1/u})
=\g_1^{-1}\cI_h(\De_i^u,\De_j^{1/u})
$$
for $i=2,4$ and $j=1,3$. Since $\g_1\ne 1$,  
$\cI_h(\De_i^u,\De_j^{1/u})=0$ for $i=2,4$ and $j=1,3$.
By (\ref{eq:transpose}),  
we have $\cI_h(\De_i^u,\De_j^{1/u})=0$ for $i=1,3$ and $j=2,4$.
To show $\cI_h(\De_i^u,\De_j^{1/u})=0$ for $(i,j)=(1,3),(2,4),(3,1),(4,2)$, 
use the map $(\rho_2)_*$. \qqed

\begin{figure}[htb] 
\centering
\includegraphics[width=8cm]{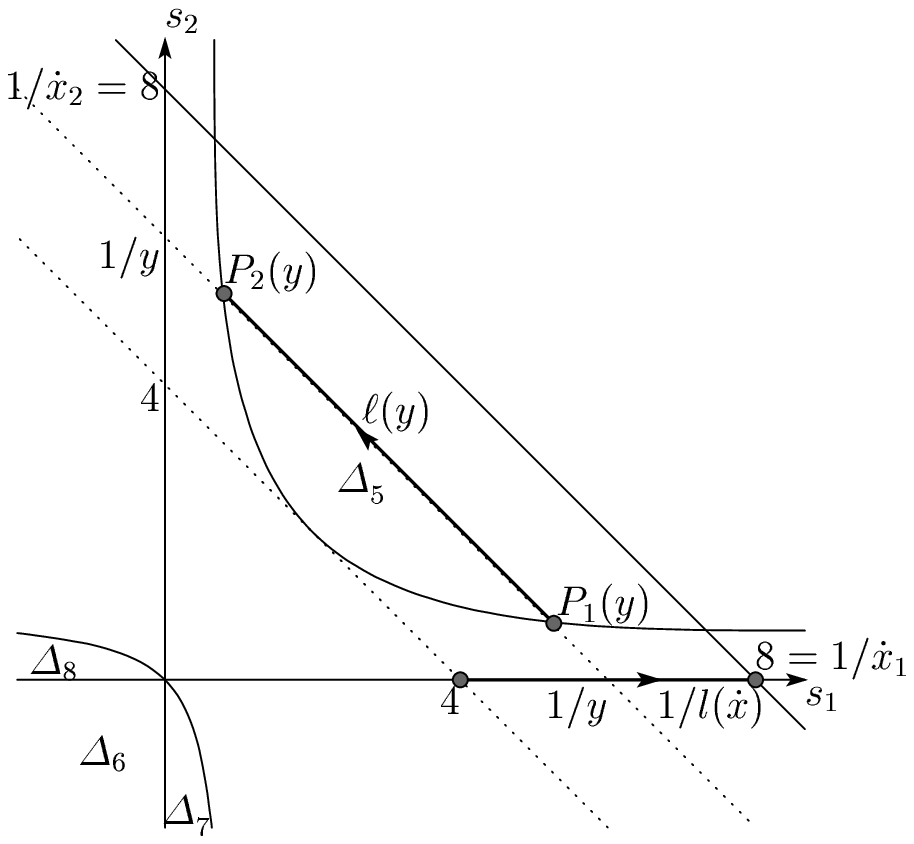}
\caption{Cycles $\De_5,\dots,\De_8$}
\label{fig:D5678}
\end{figure}

\begin{remark}
\label{rem:V1-V2}
The eigenspace $V_1^u$  of $\cM_1$ with eigenvalue $1$ is 
spanned by $\De_1^u$ and $\De_3^u$.
The eigenspace of $\cM_1$ with eigenvalue $1/\g_1$ is characterized by 
$$
\{\De^u\in H_2(\cC_\bu(\C_{\dot x}^2,u))\mid  
\cI_h(\De^u,\De_1^{1/u})=\cI_h(\De^u,\De_3^{1/u})=0\}.
$$
The eigenspace $V_2^u$  of $\cM_2$ with eigenvalue $1$ is 
spanned by $\De_1^u$ and $\De_2^u$.
The eigenspace of $\cM_2$ with eigenvalue $1/\g_2$ is characterized by 
$$
\{\De^u\in H_2(\cC_\bu(\C_{\dot x}^2,u))\mid  
\cI_h(\De^u,\De_1^{1/u})=\cI_h(\De^u,\De_2^{1/u})=0\}.
$$
Note that the linear transformation $\cM_i$ $(i=1,2)$ is determined by 
the subspace $V_i^u$, the eigenvalue $1/\g_i$ 
and the intersection form $\cI_h$, under the condition 
$c_i\notin \Z$ 
when we assign complex values to the parameters. 
\end{remark}
We characterize the linear transformation $\cM_3$ by determining 
its eigenvalues and eigenspaces. 
The following is the key lemma of this section.  
\begin{lemma}
\label{lem:M3}
We have
$$\cM_3(\De_5^u)=-\mu_3\mu_4\De_5^u
=-\dfrac{\g_1\g_2}{\a\b}\De_5^u,\quad 
\cM_3(\De^u)=\De^u
$$
for any $\De^u\in (\De_5^{1/u})^\perp
=\{\De^u\in H_2(\cC_\bu(\C_{\dot x}^2,u))\mid  
\cI_h(\De^u,\De_5^{1/u})=0\}.$
\end{lemma}
\proof
We express $\De_5$ in terms of the coordinates $s=(s_1,s_2)=(t_1/x_1,t_2/x_2)$.  
Since $L(t)$ and $Q(t,x)$ are expressed as 
$$1-{s_1}{x_1}-{s_2}{x_2},\quad x_1x_2(s_1s_2-s_1-s_2),
$$
in terms of these coordinates, we set 
$$L(s,x)=1-{s_1}{x_1}-{s_2}{x_2},\quad Q(s)=s_1s_2-s_1-s_2.$$
The intersection points $P_1$ and $P_2$ of the curves defined by 
$L(s,x)=0$ and $Q(s)=0$  are 
$$\left(\frac{1\!+\! x_1\!-\! x_2\!+\!\sqrt{R(x)}}{2x_1},
\frac{1\!-\!x_1\!+\!x_2\!-\!\sqrt{R(x)}}{2x_2}\right),\quad 
\left(\frac{1\!+\!x_1\!-\!x_2\!-\!\sqrt{R(x)}}{2x_1},
\frac{1\!-\! x_1\!+\! x_2\!+\!\sqrt{R(x)}}{2x_2}\right).$$
Note that $R(x)=1-4x_1$ for $x=(x_1,x_1)\in \rho_3$. 
When $x_1=x_2=1/4$, $R(x)$ vanishes and $Q(s)=0$ is tangent to $L(s,x)=0$. 
For $\dot x=(1/8,1/8)$, we regard $\De_5$ as 
$$\bigcup_{y\in l(\dot x_1)}\ell(y),$$ 
where 
$l(\dot x_1)$ is the segment connecting $1/4$ and $\dot x_1=1/8$, and 
$\ell(y)$ is the segment connecting the intersection points 
of $L(s,x)=0$ and $Q(s)=0$ for $x=(y,y)$ 
with $y\in l(\dot x_1)$; see Figure \ref{fig:D5678}.
For a fixed $x=(x_1,x_1)$ in the loop $\rho_3$,  
the segment $l(x_1)$ is expressed as
$$\frac{1}{4}+(x_1-\frac{1}{4})q_1$$
by a parameter $q_1\in [0,1]$. 
For an element $y=1/4+(x_1-1/4)q_1\in l(x_1)$, the segment $\ell(y)$ 
is expressed as 
$$P_1(y)+(P_2(y)-P_1(y))q_2,$$
by a parameter $q_2\in [0,1]$,
where $P_1(y)$ and $P_2(y)$ are the intersection points $P_1$ and $P_2$ 
for $x=(y,y)$.
Hence $\De_5$ is expressed by $(q_1,q_2)\in [0,1]\times[0,1]$ as
\begin{equation}
\label{eq:Delta5}
(s_1,s_2)=
\left(\frac{2(1+(1-2q_2)\sqrt{(1-4x_1)q_1})}{1-(1-4x_1)q_1},
\frac{2(1-(1-2q_2)\sqrt{(1-4x_1)q_1})}{1-(1-4x_1)q_1}\right)
\end{equation}
for a fixed $x=(x_1,x_1)$ in the loop $\rho_3$. 

By the continuation of $\sqrt{1-4x_1}$ along the loop $\rho_3$, 
its sign changes.
We regard this sign change in the deformation of $\De_5$ along $\rho_3$ 
as a bijection of $\De_5$ with the reversing orientation given by 
$$r:[0,1]\times[0,1]\ni (q_1,q_2)\mapsto (q_1,1-q_2)\in [0,1]\times[0,1].$$

We deform the pull-backs of $s_1$, $s_2$, $L(s,x)$, and $Q(s)$ 
to $[0,1]\times[0,1]$ 
by (\ref{eq:Delta5}) along $\rho_3$ and apply $r$ to them.
It is easy to see that those of $s_1$ and $s_2$ are invariant 
under the deformation and the action.
Since those of $L(s,x)$ and $Q(s)$ are expressed as 
$$
\frac{(1-q_1)(1-4x_1)}{1-(1-4x_1)q_1},\quad 
\frac{16q_1q_2(1-q_2)(1-4x_1)}{(1-q_1(1-4x_1))^2},
$$
their arguments increase by $2\pi$ under the deformation,  
and they are invariant under $r$.
Thus the pull-back of $s_1^{\l_1}s_2^{\l_2}L(s,x)^{\l_3}Q(s)^{\l_4}$ 
to $[0,1]\times[0,1]$  by (\ref{eq:Delta5}) 
is multiplied by $\mu_3\mu_4$ under the deformation along $\rho_3$ 
and the action $r$.
By considering the orientation of $\De_5$, we have 
$$\cM_3(\De_5^u)=-\mu_3\mu_4\De_5^u.$$

It is easy to see by Figure \ref{fig:D5678} that  three chambers 
\begin{eqnarray*}
\De_6&=&\{(s_1,s_2)\in \R^2\mid s_1,s_2<0\},\\
\De_7&=&\{(s_1,s_2)\in \R^2\mid s_1, Q(s)>0, s_2<0\},\\
\De_8&=&\{(s_1,s_2)\in \R^2\mid s_2, Q(s)>0,s_1<0\}
\end{eqnarray*} 
are invariant under the deformation along $\rho_3$.  
 Thus the elements $\De_i^u$ $(i=6,7,8)$ of $H_2(\cC_\bu(\C_{\dot x}^2,u))$ 
corresponding to $\De_i$ 
are eigenvectors of $\cM_3$ with eigenvalue $1$.
Since they do not intersect $\De_5$ topologically,  
they belong to $(\De_5^{1/u})^\perp$.
To show that they are linearly independent, 
we compute the intersection numbers
$$H_{ij}=\cI_h(\De_i^u,\De_j^{1/u})\quad (6\le i,j\le 8):$$
\begin{eqnarray*}
H_{66}&=&1+\frac{1}{\mu_0-1}+\frac{1}{\mu_1-1}+\frac{1}{\mu_2-1}\\
& &
+\frac{\mu_{12}-1}{(\mu_{124}-1)(\mu_1-1)(\mu_2-1)}
+\frac{\mu_{01}-1}{(\mu_{014}-1)(\mu_0-1)(\mu_1-1)}
+\frac{\mu_{02}-1}{(\mu_{024}-1)(\mu_0-1)(\mu_2-1))},\\
H_{67}&=&-\frac{1}{\mu_1-1}
\left(1+\frac{1}{\mu_{124}-1}+\frac{1}{\mu_{014}-1}\right),\\
H_{68}&=&-\frac{1}{\mu_2-1}
\left(1+\frac{1}{\mu_{124}-1}+\frac{1}{\mu_{024}-1}\right) ,\\
H_{77}&=&1+\frac{1}{\mu_1-1}+\frac{1}{\mu_4-1}
+\frac{\mu_{14}-1}{(\mu_{124}-1)(\mu_1-1)(\mu_4-1)}
+\frac{\mu_{14}-1}{(\mu_{014}-1)(\mu_1-1)(\mu_4-1)},\\
H_{78}&=&-\frac{\mu_1\mu_4}{(\mu_4-1)(\mu_{124}-1)},\\ 
H_{88}&=&1+\frac{1}{\mu_2-1}+\frac{1}{\mu_4-1}
+\frac{\mu_{24}-1}{(\mu_{124}-1)(\mu_2-1)(\mu_4-1)}
+\frac{\mu_{24}-1}{(\mu_{024}-1)(\mu_2-1)(\mu_4-1)},
\end{eqnarray*}
and $H_{ji}=H_{ij}^\vee$ for $6\le i< j\le 8$, 
where
$$\mu_0=\frac{1}{\mu_1\mu_2\mu_3\mu_4^2}=\a,\quad \mu_{ij}=\mu_i\mu_j,\quad
\mu_{ijk}=\mu_i\mu_j\mu_k.$$
Since 
$$\det(H_{ij})_{6\le i,j\le 8}=\frac{\b^2(\a-\g_1\g_2)^2(\a\b+\g_1\g_2)}
{(\a-1)(\a-\g_1)(\a-\g_2)(\b-1)^2(\b-\g_1)(\b-\g_2)(\b-\g_1\g_2)},$$
if $\a\b+\g_1\g_2\ne0$ 
when we assign complex values to the parameters, 
then they span the eigenspace of $\cM_3$ with eigenvalue $1$ and 
the space $(\De_5^{1/u})^\perp$.
\qqed

To represent $\cM_3$ by a matrix, 
we express $\De_5^u$ by a linear combination of $\De_1^u,\dots,\De_4^u$. 

\begin{lemma}
\label{lem:V-cycle}
We have 
\begin{eqnarray*}
\cI_h(\De_5^u ,\De_1^{1/u})&=&
\dfrac{1-(\mu_1\mu_4)(\mu_2\mu_4)(\mu_3)}
{(1-\mu_1\mu_4)(1-\mu_2\mu_4)(1-\mu_3)}=
\dfrac{-(1-\a)\g_1\g_2}{(\a-\g_1\g_2)(1-\g_1)(1-\g_2)},\\
\cI_h(\De_5^u ,\De_2^{1/u})&=&\cI_h(\De_5^u ,\De_3^{1/u})
=\dfrac{\mu_3\mu_4}{(1-\mu_3)(1-\mu_4)}
=\dfrac{-\g_1\g_2}{(\a-\g_1\g_2)(1-\b)},
\\
\cI_h(\De_5^u ,\De_4^{1/u})&=&
\dfrac{1-(\mu_1\mu_4)^{-1}(\mu_2\mu_4)^{-1}(\mu_4)}
{(1-(\mu_1\mu_4)^{-1})(1-(\mu_2\mu_4)^{-1})(1-\mu_4)}=
\dfrac{-(\b-\g_1\g_2)}{(1-\b)(1-\g_1)(1-\g_2)}.
\end{eqnarray*}
The twisted cycle $\De_5^u $ is expressed as
$$\De_1^u 
-\frac{\g_2(\a-\g_1)(\b-\g_1)}{\a\b(1-\g_1)(1-\g_2)}\De_2^u  
-\frac{\g_1(\a-\g_2)(\b-\g_2)}{\a\b(1-\g_1)(1-\g_2)}\De_3^u 
+\De_4^u,
$$
this leads to 
$$
\cI_h(\De_5^u ,\De_5^{1/u})
=\dfrac{1+\mu_3\mu_4}{(1-\mu_3)(1-\mu_4)}
=\frac{-(\a\b+\g_1\g_2)}{(\a-\g_1\g_2)(1-\b)}.
$$
\end{lemma}

\proof
By the results in Section 3.4 of Chapter VIII of \cite{Y2},
we can compute the intersection numbers 
$\cI_h(\De_5^u ,\De_i^{1/u})$ for $i=2,3$.  
Among the components of $\De_1$,  only $\triangle$  intersects 
with $\sqrt{-1}\R^2_x$ at $(\sqrt{x_1},\sqrt{x_2})$.  
Since their topological intersection number at this point is $-1$,
we have
$$(\sqrt{-1}\R^2_x)^{1/u}=
\dfrac{(1-\g_1)(1-\g_2)(\a-\g_1\g_2)}{(1-\a)\g_1\g_2}
\De_1^{1/u}$$
by (\ref{eq:A-rep}). 
This implies that 
$\cI_h(\De_5^u ,\De_1^{1/u})=
\dfrac{-(1-\a)\g_1\g_2}{(\a-\g_1\g_2)(1-\g_1)(1-\g_2)}$.
We can evaluate 
the intersection number $\cI_h(\De_5^u ,\De_4^{1/u})$ in a similar way.
Lemma \ref{lem:intH} together with Lemma \ref{lem:H-int} imply 
the expression of $\De_5^u $  
as a linear combination of $\De_i^u $ $(i=1,\dots,4).$
\qqed

\begin{remark}
\label{rem:V3}
\begin{enumerate}
\item 
The eigenspace of $\cM_3$ with eigenvalue $1$ is characterized by 
$\De_5^u$ and the intersection form $\cI_h$.
\item
If $\a\b+\g_1\g_2=0$, then $\cI_h(\De_5^u,\De_5^{1/u})=0$. 
In this case, the $3$-dimensional space $(\De_5^{1/u})^\perp$ 
contains the cycle $\De_5^u$ and 
coincides with the eigenspace of $\cM_3$ with eigenvalue $1$.
Since $H_2(\cC_\bu(\C_{\dot x}^2,u))$ is not spanned by 
eigenvectors of $\cM_3$, its representation is not diagonalizable.
\end{enumerate}
\end{remark}

\begin{proposition}
\label{prop:M3}
With respect to the basis $\tr(\De_1^u,\De_2^u,\De_3^u,\De_4^u)$,
$\cM_3$ is represented by the matrix
$$
\id_4-(1+\g_1\g_2\a^{-1}\b^{-1})
\frac{H \tr e_5^\vee  e_5}{e_5 H \tr e_5^\vee}
=\id_4-\frac{(\b-1)(\a-\g_1\g_2)}{\a\b}
H\tr e_5^\vee e_5,
$$
where $\id_4$ is the unit matrix of size $4$, and 
\begin{eqnarray*}
e_5&=&\left(
1,
-\dfrac{\g_2(\a-\g_1)(\b-\g_1)}{\a\b(\g_2-1)(\g_1-1)}, 
-\dfrac{\g_1(\a-\g_2)(\b-\g_2)}{\a\b(\g_2-1)(\g_1-1)}, 
1
\right),\\
e_5^\vee&=&
\left(1,
-\dfrac{(\a-\g_1)(\b-\g_1)}{\g_1(\g_1-1)(\g_2-1)},
-\dfrac{(\a-\g_2)(\b-\g_2)}{\g_2(\g_1-1)(\g_2-1)},
1
\right),
\end{eqnarray*}
corresponding to $\De_5^u$ and $\De_5^{1/u}$
by the expression in Lemma \ref{lem:V-cycle}.
\end{proposition}
\proof 
We set $M=\id_4-(1+\g_1\g_2\a^{-1}\b^{-1})H
\tr e_5^\vee  (e_5 H \tr e_5^\vee)^{-1} e_5$.
Since
$$\cI_h(\De^u,\De_5^{1/u})=(d_1,\dots,d_4)H\tr e_5^\vee,$$ 
for $\De^u=(d_1,\dots,d_4)\tr(\De_1^u,\De_2^u,\De_3^u,\De_4^u),$
we have
\begin{eqnarray*}
e_5M&=&e_5 -(1+\g_1\g_2\a^{-1}\b^{-1}) e_5H
\tr e_5^\vee  (e_5 H \tr e_5^\vee)^{-1} e_5=
-\frac{\g_1\g_2}{\a\b} e_5,\\
 (d_1,\dots,d_4)M&=&(d_1,\dots,d_4),
\end{eqnarray*}
for $(d_1,\dots,d_4)$ satisfying $(d_1,\dots,d_4)H\tr e_5^\vee=0$.
Thus the eigenvalues of $M$ are $-\g_1\g_2/(\a\b)$ and $1$,  
$e_5$ is an eigenvector with eigenvalue $-\g_1\g_2/(\a\b)$, and 
the eigenspace with eigenvalue $1$ is characterized by the equality
$(d_1,\dots,d_4)H\tr e_5^\vee=0$. 
Since $e_5$ corresponds to $\De_5$ and 
$(d_1,\dots,d_4)H\tr e_5^\vee=\cI_h(\De^u,\De_5^{1/u})$ 
for $\De^u=d_1\De_1^u+\cdots+d_4\De_4^u$, 
the linear transformation represented by $M$  coincides with 
$\cM_3$ by Lemma \ref{lem:M3}. 
Note that 
$\dfrac{1+\g_1\g_2\a^{-1}\b^{-1}}{e_5 H \tr e_5^\vee}=
\dfrac{(\b-1)(\a-\g_1\g_2)}{\a\b}$
by Lemma \ref{lem:V-cycle}. The representation matrix of $\cM_3$ 
on the right-hand side is valid even in the case $\a\b+\g_1\g_2=0$. 
\qqed

Note that $\cM_1,\cM_2$, and $\cM_3$ are represented by the matrices 
in Propositions \ref{prop:M1M2} and \ref{prop:M3} 
with respect to the basis $\tr(\De_1^u,\De_2^u,\De_3^u,\De_4^u)$. 
However, this basis degenerates when we assign an integer to $c_i$ $(i=1,2)$.
For example, if $c_1=1$, then $\g_1=1$ and $\cM_1$ is represented by 
the unit matrix; we see that this expression is not valid in this case.
Hence we give expressions of $\cM_1,\cM_2$, and $\cM_3$ in terms of 
the intersection form $\cI_h$, which are 
independent of the choice of a basis of $H_2(\cC_\bu(\C_{\dot x}^2,u))$
and are valid even for integer values of $c_1,c_2$.
As we have mentioned in Remarks \ref{rem:V1-V2} and  \ref{rem:V3}, 
$\cM_i$ are determined by the eigenspaces $V_1^u$, $V_2^u$, 
the eigenvector $\De_5^u$, and the intersection form $\cI_h$. 
We take a basis of $H_2(\cC_\bu(\C_{\dot x}^2,u))$ consisting of 
bases of these subspaces. 
We set 
$$\Hat{\De}^u_{1234}
=\tr(\Hat{\De}_1^u,\Hat{\De}_2^u,\Hat{\De}_3^u,\Hat{\De}_4^u)=
P\tr(\De_1^u,\De_2^u,\De_3^u,\De_5^u),$$
where 
$$P=\begin{pmatrix}
\dfrac{\a\b(1-\g_1)(1-\g_2)}{(1-\a)(1-\b)\g_1\g_2}&0&0&0\\
\dfrac{-\a\b(1-\g_2)}{(1-\a)(1-\b)\g_2}&\dfrac{\g_1}{1-\g_1}&0&0\\
\dfrac{-\a\b(1-\g_1)}{(1-\a)(1-\b)\g_1}&0&\dfrac{\g_2}{1-\g_2}&0\\
0&0&0&1\\
\end{pmatrix}.
$$
\begin{lemma}
\label{lem:c-integer}
The integrals 
$$\Hat f(x)=\int_{\Hat{\De}_i}u(t,x)\f_1\quad (i=1,2,3)$$
are well defined even in the case $c_1,c_2\in \Z$ when 
we assign complex values to the parameters. 
\end{lemma}
\proof
By Lemma \ref{lem:series-rep}, we have 
\begin{eqnarray*}
\Hat f_1(x)
&=&G_4\sum_{n\in \N^2} \frac{\G(a+n_1+n_2)\G(b+n_1+n_2)}
{\G(c_1+n_1)\G(c_2+n_2)\G(1+n_1)\G(1+n_2)}x_1^{n_1}x_2^{n_2},\\
f_2(x)
&=&G_4\sum_{n\in \N^2} \frac{\G(a+1-c_1+n_1+n_2)\G(b+1-c_1+n_1+n_2)}
{\G(2-c_1+n_1)\G(c_2+n_2)\G(1+n_1)\G(1+n_2)}x_1^{n_1+1-c_1}x_2^{n_2},\\
f_3(x)
&=&G_4\sum_{n\in \N^2} \frac{\G(a+1-c_2+n_1+n_2)\G(b+1-c_2+n_1+n_2)}
{\G(c_1+n_1)\G(2-c_2+n_2)\G(1+n_1)\G(1+n_2)}x_1^{n_1}x_2^{n_2+1-c_2},\\
\Hat f_2(x)
&=&G_4\frac{\g_1}{1-\g_1}(f_2(x)-\Hat f_1(x)),\\
\Hat f_3(x)
&=&G_4\frac{\g_2}{1-\g_2}(f_3(x)-\Hat f_1(x)),
\end{eqnarray*}
where $G_4=\G(1-b)\G(c_1+c_2-a-1)e^{\pi\sqrt{-1}(a+b-c_1-c_2)}$.
It is clear that $\Hat f_1(x)$ is well defined 
for $c_1,c_2\in \Z$. We claim that 
$$\lim_{c_1\to m}\frac{f_2(x)-\Hat f_1(x)}{c_1-m}$$
converges to a nonzero function for any $m\in \Z$.
Let $m$ be a fixed integer, and put $c_1=m-\e$.
Then $f_2(x)/G_4$ is 
$$\sum_{\substack{n'_1\ge 1-m \\ n_2\ge 0}} 
\frac{\G(a+n'_1+n_2+\e)\G(b+n'_1+n_2+\e)}
{\G(1+n'_1+\e)\G(c_2+n_2)\G(n'_1+m)\G(1+n_2)}x_1^{n'_1+\e}x_2^{n_2},
$$
where $n'_1=n_1+1-m$. 
If $m\ge 2$, then we have
$$\lim_{\e\to 0}\frac{1}{\G(1+n'_1+\e)}=0$$
for $1-m\le n_1'<0$. 
If $m\le 0$, then the terms $1/\G(c_1+n_1)$ ($0\le n_1\le -m$)
in the series expressing $\Hat f_1(x)$ converge to $0$ as $c_1\to m$.
Thus 
$f_2(x)$ converges to $\Hat f_1(x)$ with $c_1=m$ as $\e\to 0$.
Since the poles of the $\G$-function are simple, we have this claim.
Similarly we can show that 
$\Hat f_3(x)$ is well defined for $c_1,c_2\in \Z$. 
\qqed

The intersection matrix  
$\Hat H=\left(\cI_h(\Hat{\De}_i^u,\Hat{\De}_j^{1/u})\right)_{1\le i,j\le 4}$ 
is given by
{\small
$$\begin{pmatrix}
\frac{-\a \b (1-\g_1) (1-\g_2)}{(1-\a)(\a-\g_1 \g_2) (1-\b)^2 }&
\frac{-\a \b (1-\g_2)}{(1-\a)(\a-\g_1 \g_2) (1-\b)^2}&
\frac{-\a\b(1-\g_1)}{(1-\a)(\a-\g_1 \g_2) (1-\b)^2 }&
\frac{-\a \b}{(\a-\g_1 \g_2) (1-\b)}\\[2mm]
\frac{\a \b \g_1(1-\g_2) }{(1-\a)(\a-\g_1 \g_2) (1-\b)^2 }&
\frac{\a\b(\a \b-\g_1)\g_1(1-\g_2)}{(1-\a)(\a-\g_1)
(\a-\g_1 \g_2) (1-\b)^2(\b-\g_1) }&
\frac{\a \b \g_1}{(1-\a)(\a-\g_1 \g_2) (1-\b)^2 }&0\\[2mm]
\frac{\a\b(1-\g_1)\g_2 }{(1-\a)(\a-\g_1 \g_2)(1-\b)^2 }&
\frac{\a \b \g_2}{(1-\a)(\a-\g_1 \g_2) (\b-1)^2 }&
\frac{\a\b(\a \b-\g_2)(1-\g_1)\g_2 }{(1-\a)(\a-\g_2) 
(\a-\g_1 \g_2)(\b-1)^2(\b-\g_2) }
&0\\[2mm]
\frac{-\g_1 \g_2}{(\a-\g_1 \g_2) (1-\b)}&0&0&
\frac{-(\a \b+\g_1 \g_2)}{(\a-\g_1 \g_2)(1-\b) }\end{pmatrix},
$$
}
and its determinant is
$$
\frac{\a^3\b^3(\b-\g_1\g_2)\g_1^2\g_2^2}
{(1-\a)(\a-\g_1)(\a-\g_2)(\a-\g_1\g_2)^3(1-\b)^5(\b-\g_1)(\b-\g_2)}.
$$
Let $\Hat H_{12}$ (resp. $\Hat H_{13}$) be the submatrix of 
$\Hat H$ made by 
entries $(1,1)$, $(1,2)$, $(2,1)$, and $(2,2)$
 (resp. $(1,1)$, $(1,3)$, $(3,1)$, and $(3,3)$).

\begin{theorem}
\label{th:monodromy}
The linear transformations $\cM_i=\cM(\rho_i)$ $(i=1,2,3)$ of 
$H_2(\cC_\bu(\C_{\dot x}^2,u))$ are expressed as
\begin{eqnarray*}
\cM_1(\De^u)&=&\dfrac{1}{\g_1}\De^u
+(1-\dfrac{1}{\g_1})\left(\cI_h(\De^u,\Hat{\De}_1^{1/u}),
\cI_h(\De^u,\Hat{\De}_3^{1/u})\right)(\Hat H_{13})^{-1}
\begin{pmatrix}
\Hat{\De}_1^u\\
\Hat{\De}_3^u
\end{pmatrix},\\
\cM_2(\De^u)&=&\dfrac{1}{\g_2}\De^u
+(1-\dfrac{1}{\g_2})\left(\cI_h(\De^u,\Hat{\De}_1^{1/u}),
\cI_h(\De^u,\Hat{\De}_2^{1/u})\right)(\Hat H_{12})^{-1}
\begin{pmatrix}
\Hat{\De}_1^u\\
\Hat{\De}_2^u
\end{pmatrix},\\
\cM_3(\De^u)&=&\De^u
-(1+\dfrac{\g_1\g_2}{\a\b})
\dfrac{\cI_h(\De^u,\De_5^{1/u})}
{\cI_h(\De_5^u,\De_5^{1/u})}\De_5^u
=
\De^u-\frac{(\b-1)(\a-\g_1\g_2)}{\a\b}
\cI_h(\De^u ,\De_5^{1/u})\De_5^u.
\\
\end{eqnarray*}
\end{theorem}

\proof
By Proposition \ref{prop:M1M2} and Lemma \ref{lem:intH}, 
the eigenspace of $\cM_1$ with eigenvalue $1$ is spanned by 
$\De_1^u$ and $\De_3^u$, and that with $\g_1^{-1}$ is 
its orthogonal complement 
$$\{\De^u\in H_2(\cC_\bu(\C_{\dot x}^2,u))\mid 
\cI_h(\De^u,\De_1^{1/u})=\cI_h(\De^u,\De_3^{1/u})=0\}.$$
The elements 
$\Hat{\De}_1^u$ and $\Hat{\De}_3^u$ belong to 
the eigenspace of $\cM_1$ with eigenvalue $1$, 
and they are linearly independent. Set
$$
\cM_1'(\De^u)=\dfrac{1}{\g_1}\De^u
+(1-\dfrac{1}{\g_1})\left(\cI_h(\De^u,\Hat{\De}_1^{1/u}),
\cI_h(\De^u,\Hat{\De}_3^{1/u})\right)(\Hat H_{13})^{-1}
\begin{pmatrix}
\Hat{\De}_1^u\\
\Hat{\De}_3^u
\end{pmatrix}.
$$
We can easily check that 
$$\cM_1'(\De^u)=\left\{
\begin{array}{ccl}
\Hat{\De}_i^u &\textrm{if} & \De^u=\Hat{\De}_i^u\ (i=1,3),\\
\dfrac{1}{\g_1}\De^u &\textrm{if} & \cI_h(\De^u,\Hat{\De}_1^{1/u})
=\cI_h(\De^u,\Hat{\De}_3^{1/u})=0,
\end{array}
\right.
$$
by the property
$$
\left(\cI_h(\De^u,\Hat{\De}_1^{1/u}),
\cI_h(\De^u,\Hat{\De}_3^{1/u})\right)(\Hat H_{13})^{-1}
=\left\{
\begin{array}{ccl}
(1,0) &\textrm{if}& \De^u=\Hat{\De}_1^{u},\\
(0,1) &\textrm{if}& \De^u=\Hat{\De}_3^{u},\\
(0,0) &\textrm{if}& \cI_h(\De^u,\Hat{\De}_1^{1/u})
=\cI_h(\De^u,\Hat{\De}_3^{1/u})=0.
\end{array}
\right.
$$
Since the eigenvalues and eigenspaces of $\cM_1$ coincide
with those of $\cM_1'$, we have $\cM_1=\cM_1'$. 
We obtain the expression of $\cM_2$ in a similar way.  
Set 
$$\cM'_3(\De^u)=\De^u
-(1+\dfrac{\g_1\g_2}{\a\b})
\frac{\cI_h(\De^u,\De_5^{1/u})}{\cI_h(\De_5^u,\De_5^{1/u})}\De_5^u.
$$
By the property
$$\frac{\cI_h(\De^u,\De_5^{1/u})}{\cI_h(\De_5^u,\De_5^{1/u})}
=\left\{
\begin{array}{ccl}
1 &\textrm{if}& \De^u=\De_5^{u},\\
0 &\textrm{if}& \De^u\in (\De_5^{1/u})^\perp,
\end{array}
\right.
$$
we see that 
$$\cM_3'(\De^u)=\left\{
\begin{array}{ccl}
-\dfrac{\g_1\g_2}{\a\b}\De_5^u &\textrm{if} & \De^u=\De_5^u,\\
\De^u &\textrm{if} & \De^u\in (\De_5^{1/u})^\perp,
\end{array}
\right.
$$
which shows $\cM_3=\cM_3'$ by Lemma \ref{lem:M3}.
The second  expression of $\cM_3$ is obtained by the equality
$\cI_h(\De_5^u ,\De_5^{1/u})
=\dfrac{-(\a\b+\g_1\g_2)}{(\a-\g_1\g_2)(1-\b)}$ 
in Lemma \ref{lem:V-cycle}.
\qqed

\begin{remark}
\begin{enumerate}
\item 
We note that when we assign integers to $c_1$ and $c_2$, although $\De_1^u$, $\De_2^u$, and $\De_3^u$ are linearly dependent, 
$\Hat{\De}_1^u$, $\Hat{\De}_2^u$, and $\Hat{\De}_3^u$
remain linearly independent.

\item
Since we have 
\begin{eqnarray*}
(\Hat H_{12})^{-1}&=&
\frac{(\a-\g_1\g_2)(1-\b)}{\a \b\g_1^2(1-\g_2)}
\begin{pmatrix}
(\a \b-\g_1)\g_1& 
(\a-\g_1)(\b-\g_1)\\
-(\a-\g_1)(\b-\g_1)\g_1& 
-(\a-\g_1)(\b-\g_1)(1-\g_1)\end{pmatrix},\\
(\Hat H_{13})^{-1}&=&\frac{(\a-\g_1 \g_2)(1-\b)}{\a\b(1-\g_1)\g_2^2}
\begin{pmatrix}
(\a \b-\g_2) \g_2& 
(\a-\g_2)(\b-\g_2)\\
-(\a-\g_2)(\b-\g_2)\g_2& 
-(\a-\g_2)(\b-\g_2)(1-\g_2)
\end{pmatrix},\\
\end{eqnarray*}
the factors $1-\g_1$ and $1-\g_2$ 
are canceled in the expression of $\cM_1$ and $\cM_2$. 
Theorem \ref{th:monodromy} is valid even in the case 
$c_1,c_2,a+b-c_1-c_2-\dfrac{1}{2}\in \Z$  
when we assign complex values to the parameters.
\end{enumerate}
\end{remark}

\begin{cor}
\label{cor:mat-rep} The linear transformations $\cM_i$ $(i=1,2,3)$ are
represented by matrices $M_i$ with respect to the basis 
$\Hat{\De}^u_{1234}=\tr(\Hat{\De}_1^u,\dots,\Hat{\De}_4^u)$ as
$\cM_i(\Hat{\De}^u_{1234})=M_i\Hat{\De}^u_{1234},$  
where 
\begin{eqnarray*}
M_1&=&\frac{1}{\g_1}\id_4+(1-\frac{1}{\g_1})
\Hat H(\tr e_1,\tr e_3)(\Hat H_{13})^{-1}
\begin{pmatrix}
e_1\\ e_3
\end{pmatrix}
=
\begin{pmatrix}
1&0&0&0\\
1&\frac{1}{\g_1}&0&0\\
0&0&1&0\\
\frac{\a \b-\g_2}{\a \b}&0&
\frac{(\a-\g_2)(\b-\g_2) }{\a \b\g_2 }&\frac{1}{\g_1}\end{pmatrix},
\\
M_2&=&
\frac{1}{\g_2}\id_4+(1-\frac{1}{\g_2})
\Hat H(\tr e_1,\tr e_2)(\Hat H_{12})^{-1}
\begin{pmatrix}
e_1\\ e_2
\end{pmatrix}
=
\begin{pmatrix}
1&0&0&0\\
0&1&0&0\\
1&0&\frac{1}{\g_2}&0\\
\frac{\a \b-\g_1}{\a \b}&\frac{(\a-\g_1) (\b-\g_1)}{\a \b\g_1 }&0&
\frac{1}{\g_2}\end{pmatrix},\\
M_3&=&\id_4 -(1+\frac{\g_1\g_2}{\a\b})
\frac{\Hat H\tr e_4e_4}{e_4\Hat H\tr e_4}
=
\begin{pmatrix}
1&0&0&-1\\
0&1&0&0\\
0&0&1&0\\
0&0&0&\frac{-\g_1 \g_2}{\a \b}\end{pmatrix},
\end{eqnarray*}
and $e_i$ is the $i$-th unit row vector of $\Z^4$.
\end{cor}

\proof
The matrix $M_3$ is obtained in the same way as in the proof of 
Proposition \ref{prop:M3}. By the expression of $\cM_1$ 
in Theorem \ref{th:monodromy}, 
we give its representation matrix with respect to the basis 
$\Hat{\De}^u_{1234}$.
Set ${\De}^u=(d_1,\dots,d_4)\Hat{\De}^u_{1234}$. 
Since 
$$
\cI_h(\De^u,\Hat{\De}_i^{1/u})=(d_1,\dots,d_4)\Hat{H}\tr e_i\quad (i=1,\dots,4)
$$
we have 
$$\left(\cI_h(\De^u,\Hat{\De}_1^{1/u}),
\cI_h(\De^u,\Hat{\De}_3^{1/u})\right)=
(d_1,\dots,d_4)\Hat{H}(\tr e_1,\tr e_3).$$
Note that 
$$\begin{pmatrix}
\Hat{\De}_1^u\\
\Hat{\De}_3^u
\end{pmatrix}=
\begin{pmatrix}
e_1\\
e_3
\end{pmatrix}\Hat{\De}^u_{1234}.
$$
Thus we have 
$$
\left(\cI_h(\De^u,\Hat{\De}_1^{1/u}),
\cI_h(\De^u,\Hat{\De}_3^{1/u})\right)(\Hat H_{13})^{-1}
\begin{pmatrix}
\Hat{\De}_1^u\\
\Hat{\De}_3^u
\end{pmatrix}
=(d_1,\dots,d_4)\Hat{H}(\tr e_1,\tr e_3)(\Hat H_{13})^{-1}
\begin{pmatrix}
e_1\\
e_3
\end{pmatrix}\Hat{\De}^u_{1234},
$$
which implies that $M_1$ is the representation matrix of $\cM_1$. 
We obtain the matrix $M_2$ in a similar way.
\qqed

\begin{remark}
\label{rem:mat-rep}
With respect to the basis $P'\tr(\De_1^u,\De_2^u,\De_3^u,\De_4^u)$ of 
$H_2(\cC_\bu(\C_{\dot x}^2,u))$ for  
$$P'=
\left( \begin {array}{cccc} 
\dfrac{ \a \b(1\!-\! \g_1 )(1\!-\! \g_2)}{(1 \!-\! \a)(1\!-\! \b)\g_1 \g_2} 
&0&0&0\\ 
\dfrac {-\a \b(1\!-\! \g_2)}{(1 \!-\! \a)(1\!-\! \b)\g_2}
&\dfrac {\g_1}{1\!-\! \g_1}
&0&0\\ 
\dfrac {-\a\b(1\!-\! \g_1)}{(1\!-\! \a )(1\!-\! \b)\g_1}&0&
\dfrac {\g_2}{ 1\!-\! \g_2}&0\\ 
\dfrac {\a \b}{ (1\!-\! \a  )(1\!-\! \b)}
&\dfrac {-\g_1\g_2}{(1\!-\! \g_1 ) (1\!-\! \g_2)}
&\dfrac {- \g_1\g_2}{(1\!-\!  \g_1)(1\!-\! \g_2)}
&\dfrac {\a \b \g_1\g_2}{ (\a \!-\! \g_1\g_2)(\b\!-\! \g_1\g_2)}
\end {array} \right),
$$
$\cM_1$, $\cM_2$, and $\cM_3$ are represented by matrices 
$$\begin{pmatrix}
1 &0 &0 &0\\
1 &\g_1^{-1} &0 &0 \\
0 &0  & 1 & 0\\
0 &0  & 1 & \g_1^{-1}
\end{pmatrix},
\quad 
\begin{pmatrix}
1 &0 &0 &0\\
0 &1 &0 &0\\
1 &0&\g_2^{-1}&0 \\
0 &1&0 & \g_2^{-1}
\end{pmatrix},
$$
$$\begin{pmatrix}
-\frac{\g_1\g_2}{\a\b}& \frac{\g_1\g_2}{\a\b}-\frac{1}{\g_1}
&\frac{\g_1\g_2}{\a\b}-\frac{1}{\g_2}
& -\frac{(\a-\g_1\g_2)(\b-\g_1\g_2)}{\a\b\g_1\g_2}
\\
0 &1 &0 &0\\
0 &0 &1 &0\\
0 &0 &0 &1
\end{pmatrix},
$$
respectively.
These representations of $\cM_i$ are also valid even in the case 
$c_1,c_2,a+b-c_1-c_2-\dfrac{1}{2}\in \Z$  
when we assign complex values to the parameters. 
\end{remark}

\section{Twisted cohomology group}
Recall that
\begin{eqnarray*}
\l_1&=&b+1-c_1,\quad \l_2=b+1-c_2,\quad 
\l_3=-a+c_1+c_2-1,
\quad \l_4=-b,\\
\frak X&=&\big\{(t,x)\in \C^2\times X\big|t_1t_2L(t)Q(t,x)\ne0\big\}
\subset (\P^1\times \P^1)\times \P^{2},\\ 
\C_x^2&=& \pr^{-1}(x),\quad \pr:\frak X\ni (t,x)\mapsto x\in X.
\end{eqnarray*}

In this section, we regard vector spaces as defined over 
the rational function field $\C(\l)=\C(\l_1,\dots,\l_4)=\C(a,b,c_1,c_2)$.
We denote the vector space of rational functions on $\P^2$ 
with poles only along $S$ by $\cO_X(*S)$. Note 
that $\cO_X(*S)$ admits the structure of an algebra over $\C(\l)$.
We set 
$$\frak S=(\P^1\times \P^1)\times \P^{2}-\frak X.$$
Let $\W^k_{\frak X}(*\frak S)$ 
be the vector space of rational $k$-forms on $\frak X$ 
with poles only along $\frak S$ and 
$\W^{p,q}_{\frak X}(*\frak S)$ be the subspace of 
$\W^{p+q}_{\frak X}(*\frak S)$
consisting of elements that are $p$-forms with respect to the variables 
$t_1,t_2$. 
We set
$$
\w=d_t\log(u(t,x))=\l_1\frac{dt_1}{t_1}+\l_2\frac{dt_2}{t_2}+
\l_3\frac{d_t L(t)}{L(t)}
+\l_4\frac{d_tQ(t,x)}{Q(t,x)}
\in \W^{1,0}_{\frak X}(*\frak S),
$$
where 
$d_t$ is the exterior derivative with respect to the variables
$t_1,t_2$. 
Note that 
$$d_tL(t)=-dt_1-dt_2,\quad d_t Q(t,x)=(t_2-x_2)dt_1+(t_1-x_1)dt_2.$$
By a twisted exterior derivative $\na=d_t+\w\wedge$ on $\frak X$,
we define quotient spaces
$$
\cH^k(\na)=\ker\left(\na:\W_{\frak X}^{k,0}(*\frak S)\to 
\W_{\frak X}^{k+1,0}(*\frak S)\right)
\big/\na\left(\W_{\frak X}^{k-1,0}(*\frak S)\right)\quad (k=0,1,2),$$
where we regard $\W_{\frak X}^{-1,0}(*\frak S)$ as the zero vector space.
Each of them  admits the structure of a vector bundle  over $X$.

We consider the structure of the fiber of $\cH^k(\na)$ at $x$.
Let $\W_{\C_x^2}^p(*x)$ be the vector space of rational $p$-forms on $\C_x^2$
with poles only along 
the pole divisor of the pull-back $\w_x=\imath_x^*(\w)$ of $\w$ 
by the map $\imath_x:\C_x^2\to \frak X$.
There is a natural map from each fiber of $\cH^k(\na)$ at $x$ to 
the rational twisted cohomology group
$$H^k(\W_{\C_x^2}^\bu(*x),\na_x)=
\ker\left(\na_x:\W_{\C_x^2}^k(*x)\to \W_{\C_x^2}^{k+1}(*x)\right)
/\na_x\left(\W_{\C_x^2}^{k-1}(*x)\right)$$
on $\C_x^2$ with respect to the twisted exterior derivative 
$\na_x=d_t+\w_x\wedge$.

\begin{fact}[\cite{AoKi},\cite{Cho}]
\label{fact:TCH}
\begin{itemize}
\item[$\mathrm{(i)}$] We have 
$$\dim H^k(\W_{\C_x^2}^\bu(*x),\na_x)=\left\{
\begin{array}{ccl}
4 &\textrm{if}&k=2,\\
0 &\textrm{if}&k=0,1.
\end{array}
\right.
$$
\item[$\mathrm{(ii)}$]
There is a canonical isomorphism  
$$\jmath_x: H^2(\W_{\C_x^2}^\bu(*x),\na_x) 
\to
H^2(\cE_c^\bu(x),\na_x)=\ker(\na_x:\cE_c^{2}(x)\to \cE_c^{3}(x))
/\na_x(\cE_c^{1}(x)),$$
where $\cE_c^k(x)$ is the vector space of smooth $k$-forms 
with compact support in $\C_x^2$.
\end{itemize}
\end{fact}

We have a twisted exterior derivative $\na^\vee=d_t-\w\wedge$
for $-\w$ and 
$$
\cH^2(\na^\vee)=
\W_{\frak X}^{2,0}(*\frak S)\big/\na^\vee(\W_{\frak X}^{1,0}(*\frak S)),\quad 
H^2(\W_{\C_x^2}^\bu(*x),\na_x^\vee)=
\W_{\C_x^2}^2(*x)/\na_x^\vee(\W_{\C_x^2}^{1}(*x)).
$$
The $\cO_X(*S)$-module $\cH^2(\na^\vee)$ 
can be regarded as a vector bundle over $X$.

For any fixed $x\in X$, we define the intersection form between 
$H^2(\W_{\C_x^2}^\bu(*x),\na)$
and 
$H^2(\W_{\C_x^2}^\bu(*x),\na^\vee)$
by 
$$\cI_c(\f_x,\f_x')=\int_{\C_x^2}\jmath_x(\f_x)\wedge \f_x'\in \C(\a),$$
where $\f_x,\f_x'\in \W_{\C_x^2}^{2}(*x)$,  $\jmath_x$ is given in 
Fact \ref{fact:TCH}. 
This integral converges since $\jmath_x(\f_x)$ is a smooth $2$-from 
on $\C_x^2$ with compact support. 
It is bilinear over $\C(\a)$.

We take four elements 
$$
\begin{array}{ll}
\f_1=d_t\log\big(\dfrac{t_1}{L(t)}\big)\wedge 
d_t\log\big(\dfrac{t_2}{L(t)}\big)
=\dfrac{dt_1\wedge dt_2}{t_1t_2L(t)},
&
\f_2=d_t\log(t_2)\wedge d_t\log(L(t))=\dfrac{dt_1\wedge dt_2}{t_2L(t)},
\\
\f_3=-d_t\log(t_1)\wedge d_t\log(L(t))=\dfrac{dt_1\wedge dt_2}{t_1L(t)},
& 
\f_4=\dfrac{t_1\wedge t_2}{L(t)Q(t,x)}
\end{array}
$$
of $\cH^2(\na)$, and we denote  
$\imath_x^*(\f_i)\in H^2(\W_{\C_x^2}^\bu(*x),\na_x)$
by $\f_{x,i}$.
Since $\na^\vee(\f_i)=0$, $\na_x^\vee(\f_{x,i})=0$, 
we can regard $\f_i$ and $\f_{x,i}$ as 
elements of 
$\cH^2(\na^\vee)$ and $H^2(\W_{\C_x^2}^\bu(*x),\na_x^\vee)$,
respectively.
The intersection numbers $\cI_c(\f_{x,i},\f_{x,j})$ 
($\f_{x,i}\in H^2(\W_{\C_x^2}^\bu(*x),\na_x)$,
$\f_{x,j}\in H^2(\W_{\C_x^2}^\bu(*x),\na_x^\vee)$ 
$1\le i,j\le 4$)  
are evaluated as follows.

\begin{theorem}
\label{th:int-c}
The intersection matrix $\big(\cI_c(\f_{x,i},\f_{x,j})\big)_{1\le i,j\le 4}$ 
is $(2\pi\sqrt{-1})^2 C$, where $C$ is a symmetric matrix with 
entries
\begin{eqnarray*}
C_{11}&=& \left(\frac{1}{\l_1}+\frac{1}{\l_2}\right)
\left(\frac{1}{\l_3}+\frac{1}{\l_{124}}\right)
=\dfrac{(-a+1+b)(2b-c_1-c_2+2)}
{(-a+c_1+c_2-1)(b-c_1+1)(b-c_2+1)(b-c_1-c_2+2)},\\
C_{12}&=&\frac{1}{\l_2\l_3}=\frac{1}{(b-c_2+1)(-a+c_1+c_2-1)},\\
C_{13}&=&\frac{1}{\l_1\l_3}=\frac{1}{(b-c_1+1)(-a+c_1+c_2-1)},\\
C_{14}&=&0,\\
C_{22}&=&\left(\frac{1}{\l_0}+\frac{1}{\l_2}\right)
\left(\frac{1}{\l_3}+\frac{1}{\l_{134}^-}\right)
=\dfrac{(c_1-1)(a+b-c_2)}{(a-1)(a-c_2)(b-c_2+1)(-a+c_1+c_2-1)},\\
C_{23}&=&\frac{-1}{\l_0\l_3}=\frac{-1}{(a-1)(-a+c_1+c_2-1)},\\
C_{24}&=&0,\\
C_{33}&=&\left(\frac{1}{\l_0}+\frac{1}{\l_1}\right)
\left(\frac{1}{\l_3}+\frac{1}{\l_{234}^-}\right)
=\dfrac{(c_2-1)(a+b-c_1)}{(a-1)(a-c_1)(b-c_1+1)(-a+c_1+c_2-1)},\\
C_{34}&=&0,\\
C_{44}&=&\frac{2}{\l_3\l_4R(x)}=\frac{2}{(-a+c_1+c_2-1)(-b)R(x)},
\end{eqnarray*}
where 
$$\begin{array}{ll}
\hspace{3mm}\l_0=-\l_1-\l_2-\l_3-2\l_4=a-1,&
\l_{124}=\l_1+\l_2+\l_4=b-c_1-c_2+2,\\
\l_{134}^-=-\l_1-\l_3-\l_4=a-c_2,&
\l_{234}^-=-\l_2-\l_3-\l_4=a-c_1.
\end{array}
$$
The determinant of $C$ is
$$
\frac{-4b}{(a-1)(a-c_1)(a-c_2)(-a+c_1+c_2-1)^3(b-c_1+1)(b-c_2+1)
(b-c_1-c_2+2)R(x)}.
$$
\end{theorem}

\proof 
We blow up $\P^1\times \P^1(\supset \C_x^2)$ at the two points 
$(0,0)$ and $(\infty,\infty)$ so that 
the pole divisor of $\w_x$ is normally crossing.
We tabulate the residue of $\w_x$ at  
each component of the pole divisor in Table \ref{tab:res}, 
where $E_0$ and $E_\infty$ are the exceptional divisors 
corresponding to the points $(0,0)$ and  $(\infty,\infty)$, respectively.
\begin{table}
\begin{center}
\begin{tabular}[hbt]{|l|l|}
\hline
component  &  residue\\
\hline
$E_{\infty}$ & $\l_{0}= a-1$       \\
\hline
$t_1=0$   &      $\l_{1}=b-c_1+1$ \\
\hline
$t_2=0$   &     $\l_{2}=b-c_2+1$ \\
\hline
$L(t)=0$  &     $\l_{3}=-a+c_1+c_2-1$ \\
\hline
$Q(t,x)=0$ &    $\l_{4}=-b$\\
\hline
$E_0$ & $\l_{124}=b-c_1-c_2+2$\\
\hline
$t_1=\infty$   & $\l_{134}^-=a-c_2$\\
\hline
$t_2=\infty$    & $\l_{234}^-=a-c_1$\\
\hline
\end{tabular}
\end{center}
\label{tab:res}
\caption{Residues of $\w_x$}
\end{table}
\begin{figure}[htb]
  \centering
\includegraphics[width=10cm]{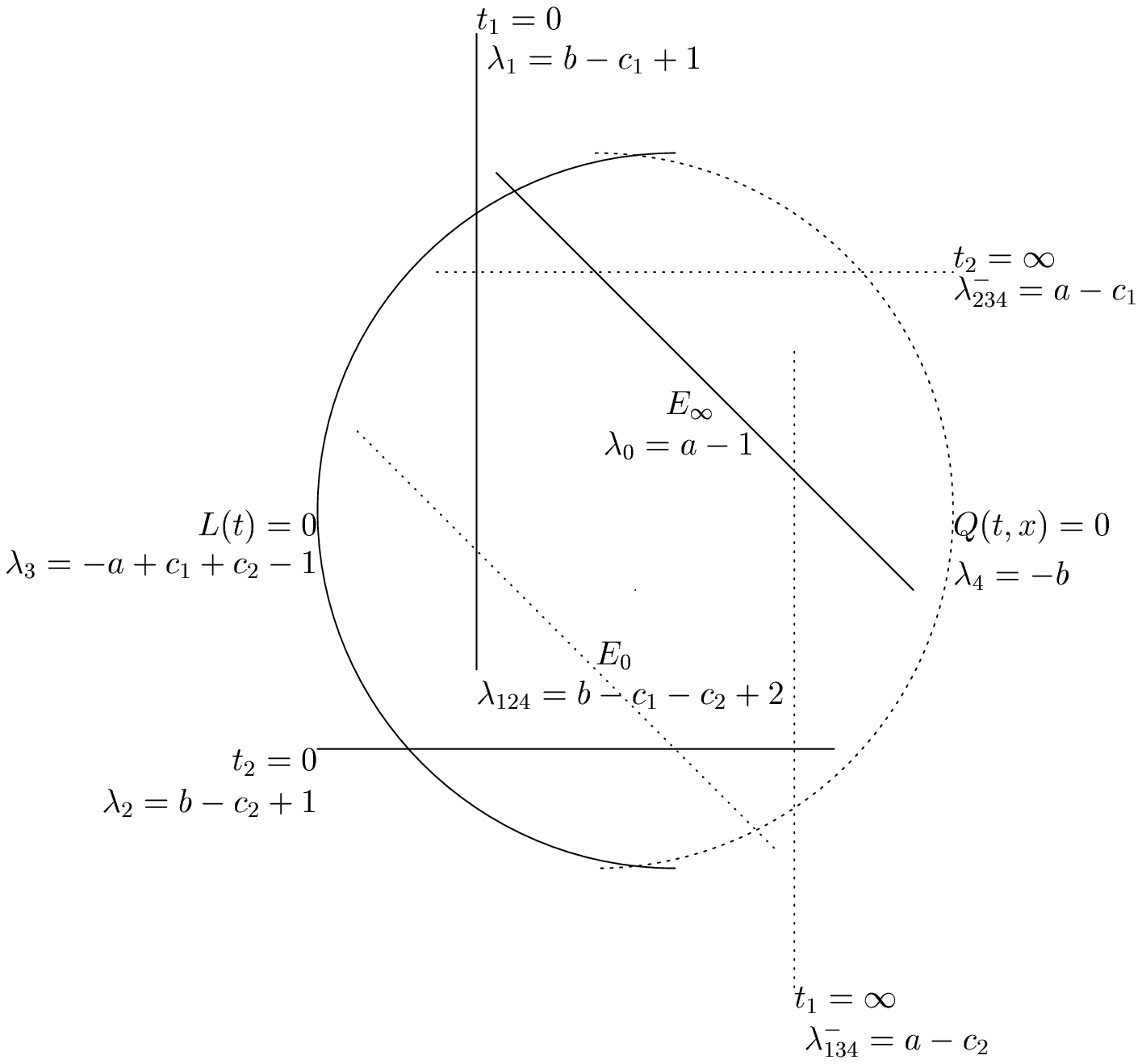}
  \caption{Pole divisor of $\w_x$}
\label{fig:exponents}
\end{figure}
To evaluate $C_{11}$, we find the intersection points 
of components of the pole divisor of $\f_{x,1}$. There are four 
points
$$\{t_1=0\}\cap E_0,\quad \{t_2=0\}\cap E_0,\quad 
\{t_1=0\}\cap \{L(t)=0\},\quad \{t_2=0\}\cap \{L(t)=0\};$$
see Figure \ref{fig:exponents}.
For every intersection point, we compute 
the reciprocal of the product of the residues of $\w_x$ along 
the components passing it. The results in Section 5 of \cite{M1}
imply that  $C_{11}$ is given by  their sum: 
$$\frac{1}{\l_1\l_{124}}+\frac{1}{\l_2\l_{124}}+
\frac{1}{\l_1\l_{3}}+\frac{1}{\l_2\l_{3}}.$$
Similarly, we can evaluate 
$C_{22}$ and $C_{33}$. 

Let us evaluate $C_{12}$. 
The intersection points of the components of the pole divisor of $\f_{x,2}$ 
are 
$$\{t_2=0\}\cap \{L(t)=0\},\quad \{t_2=0\}\cap \{t_1=\infty\},\quad 
\{L(t)=0\}\cap E_\infty,\quad \{t_1=\infty\}\cap E_\infty;$$
$\{t_2=0\}\cap \{L(t)=0\}$ is the common intersection point
of the pole divisors of $\f_{x,1}$ and $\f_{x,2}$ .
By regarding $L(t)$ and $t_2$ as local coordinates around this point, 
we express $\f_{x,1}$ and $\f_{x,2}$ in terms of them:
$$\f_{x,1}=-\frac{dL(t)\wedge dt_2}{(1-L(t)-t_2)t_2L(t)},\quad 
\f_{x,2}=-\frac{dL(t)\wedge dt_2}{t_2L(t)}.$$
Since $1/(1-L(t)-t_2)=1$ for $(L(t),t_2)=(0,0)$, 
the intersection number $C_{12}$ is given by 
the reciprocal of the product of the residues of $\w_x$ along 
the components passing the point $(L(t),t_2)=(0,0)$, that is 
$1/(\l_2\l_3)$.
Similarly, we can evaluate $C_{13}$. 
To evaluate $C_{23}$, we express $\f_{x,2}$ and  $\f_{x,3}$ in terms of 
coordinates $s_1=1/t_1$, $s_2=t_2/t_1$ around 
$\{L(t)=0\}\cap E_\infty$ represented by $(s_1,s_2)=(0,-1)$. 
Since 
$$\f_{x,2}=\frac{-ds_1\wedge ds_2}{s_1(s_1-1-s_2)},\quad 
\f_{x,3}=\frac{-ds_1\wedge ds_2}{s_1s_2(s_1-1-s_2)},\quad
\left[s_2\right]_{(s_1,s_2)=(0,-1)}=-1,$$
and the residue of $\w_x$ along $\{L(t)=0\}$ and 
that along $E_\infty$ are $\l_3$ and $\l_0$, respectively, 
we have $C_{23}=-1/(\l_0\l_3)$. 

The pole divisor of $\f_4$ consists of $L(t)=0$ and $Q(t,x)=0$. 
They intersect at the two points $P_1$ and $P_2$.
Since the pole divisor of $\f_{x,i}$ $(i=1,2,3)$ does not contain  
$Q(t,x)=0$, we have $C_{i4}=0$ for $i=1,2,3$. 
To compute $C_{44}$, we express $\f_4$ around 
the intersection points $P_1$ and $P_2$ 
in terms of the local coordinates $L(t)$ and $Q(t,x)$.
A straightforward calculation implies 
$$\f_4=\frac{(-1)^idL(t)\wedge dQ(t,x)}{L(t)Q(t,x)\sqrt{R(x)+L(t)^2
-2(1-x_1-x_2)L(t)-4Q(t,x)}}$$
around $P_i$ ($i=1,2$), 
where the function ${(-1)^i}/{\sqrt{R(x)+L(t)^2-2(1-x_1-x_2)L(t)-4Q(t,x)}}$ 
is a single-valued holomorphic function around $P_i$ 
with value ${(-1)^i}/{\sqrt{R(x)}}$ at this point. 
We have 
$$C_{44}=\frac{1}{\l_3\l_4}\frac{{-1}}{{\sqrt{R(x)}}}\frac{{-1}}{{\sqrt{R(x)}}}
+\frac{1}{\l_3\l_4}\frac{{1}}{{\sqrt{R(x)}}}\frac{{1}}{{\sqrt{R(x)}}}
=\frac{2}{\l_3\l_4R(x)}.$$

The determinant of $C$ is obtained by a straightforward calculation.
\qqed

Note that the matrix $C$ is well defined, and 
$\det(C)\ne0$ for any $x\in X$ under our assumption.
The natural map from each fiber of $\cH^2(\na)$ at $x$ to 
$H^2(\W_{\C_x^2}^\bu(*x),\na_x)$ is surjective.
The $\C(\l)$-span of the classes of $\f_1,\dots,\f_4
\in \cH^2(\na)$ (resp. $\in \cH^2(\na^\vee)$) is denoted by  
$\cH^2_{\C(\l)}(\na)$  (resp. $\cH^2_{\C(\l)}(\na^\vee)$). 
The intersection form $\cI_c$ is regarded as 
a map from $\cH^2_{\C(\l)}(\na)\times \cH^2_{\C(\l)}(\na^\vee)$ 
to $\cO(*S)$.

\section{Twisted period relations}
Among $H_2(\cC_\bu(\C_{x}^2,u))$, 
$H_2(\cC_\bu(\C_{x}^2,1/u))$,
$H^2(\W_{\C_x^2}^\bu(*x),\na_x)$, and 
$H^2(\W_{\C_x^2}^\bu(*x),\na_x^\vee)$, 
there are the intersection pairings $\cI_h$ and $\cI_c$, and 
the pairings which yield solutions of $\cF_4$ with various parameters.
We have two isomorphisms from $H_2(\cC_\bu(\C_{x}^2,u))$ 
to $H^2(\W_{\C_x^2}^\bu(*x),\na_x^\vee)$ by regarding 
them as the dual spaces of $H_2(\cC_\bu(\C_{x}^2,1/u))$
and those of $H^2(\W_{\C_x^2}^\bu(*x),\na_x)$.
As is shown in \cite{KY}, these isomorphisms coincide.
This compatibility implies the following. 
\begin{theorem}
\label{th:TPR}
The intersection matrices $H$ and $(2\pi\sqrt{-1})^2C$ and 
the period matrices 
$$\varPi_+(x)= \left(\int_{\De_j} u\f_{x,i}\right)_{0\le i,j\le 4},
\quad 
\varPi_-(x)=\left(\int_{\De_j}\! (1/u)\f_{x,i}\right)_{0\le i,j\le 4}
$$
satisfy 
\begin{equation}
\label{eq:TPR}
\varPi_+(x)\tr H^{-1}\tr \varPi_-(x)=(2\pi\sqrt{-1})^2C.
\end{equation}
\end{theorem}

\begin{cor}
\label{cor:TPR}
The identity 
(\ref{eq:TPR}) implies twisted period relations 
\begin{eqnarray*}
& &\frac{1-a}{1-a_{12}} 
F_4(a,b,c_1,c_2;x)F_4(2-a, -b, 2-c_1, 2-c_2;x)\\
&-&\frac{b(1-a_1)}{b_1(1-a_{12})}
F_4(a_1,b_1,2-c_1,c_2;x)F_4(2-a_1,-b_1,c_1,2-c_2;x)\\
&-&\frac{b(1-a_2)}{b_2(1-a_{12})}
F_4(a_2,b_2,c_1,2-c_2;x)F_4(2-a_2,-b_2,2-c_1,c_2;x)\\
&+&\frac{b}{b_{12}}
F_4(a_{12},b_{12},2-c_1,2-c_2;x)F_4(2-a_{12},-b_{12},c_1,c_2;x)
\\
&=&
\frac{(1-a+ b)(b_1+b_2)(1-c_1)(1-c_2)}
{(1-a_{12})b_1b_2b_{12}},
\end{eqnarray*}
\begin{eqnarray*}
& &\frac{1-a}{1-a_{12}}
F_4(a,b+1,c_1,c_2;x)F_4(2-a,1-b,2-c_1,2-c_2;x)\\
&-&\frac{b_1(1-a_1)}{b(1-a_{12})}
F_4(a_1,b_1+1,2-c_1,c_2;x)F_4(2-a_1,1-b_1,c_1,2-c_2;x)\\
&-&
\frac{b_2(1-a_2)}{b(1-a_{12})}
F_4(a_2,b_2+1,c_1,2-c_2;x)F_4(2-a_2,1-b_2,2-c_1,c_2;x)\\
&+&
\frac{b_{12}}{b}
F_4(a_{12},b_{12}+1,2-c_1,2-c_2;x)
F_4(2-a_{12},1-b_{12},c_1,c_2;x)\\
&=&\frac{2(1-c_1)(1-c_2)}{(1-a_{12})(-b)R(x)},
\end{eqnarray*}
\begin{eqnarray*}
& &\frac{1-a}{1-a_{12}}F_4(a, b, c_1, c_2;x)F_4(2-a, 1-b, 2-c_1, 2-c_2;x)\\
&-&\frac{1-a_1}{1-a_{12}}
F_4(a_1,b_1,  2-c_1, c_2;x)F_4(2-a_1,1-b_1,  c_1, 2-c_2;x)\\
&-&\frac{1-a_2}{1-a_{12}}
F_4(a_2, b_2, c_1, 2-c_2;x)F_4(2-a_2,1-b_2, 2-c_1, c_2;x)\\
&+&
F_4(a_{12},b_{12},2-c_1,2-c_2;x)F_4(2-a_{12},1-b_{12},c_1,c_2;x)\\
&=&0,
\end{eqnarray*}
where
$$\begin{array}{lll}
a_1=a-c_1+1,& a_2=a-c_2+1,& a_{12}=a-c_1-c_2+2,\\
b_1=b-c_1+1,& b_2=b-c_2+1,& b_{12}=b-c_1-c_2+2.
\end{array}
$$

\end{cor}
\proof
Compare the $(1,1)$-entries of the both sides of (\ref{eq:TPR}). 
Then we have 
\begin{equation}
\label{eq:TPR11}
(f_1(x),\dots,f_4(x)) \tr H^{-1} \tr(f_1^\vee(x),\dots,f_4^\vee(x))
=\cI_c(\f_{x,1},\f_{x,1}),
\end{equation}
where
\begin{eqnarray*}
f_1^\vee(x)&=&\frac{\G(c_1-1)\G(c_2-1)\G(1-c_1-c_2+a)}{\G(-1+a)} 
F_4(2-a, -b, 2-c_1, 2-c_2;x),\\
f_2^\vee(x)&=&
\frac{\G(c_1-b-1)\G(c_1-a+1)\G(1+b)\G(1-c_1-c_2+a)}{\G(c_1)\G(2-c_2)}\\
& &\times e^{-\pi\sqrt{-1}(a+b-c_1-c_2)}x_1^{c_1-1}
F_4(c_1-a+1, c_1-b-1, c_1, 2-c_2;x),\\
f_3^\vee(x)&=&
\frac{\G(c_2-a+1)\G(c_2-b-1)\G(1+b)\G(1-c_1-c_2+a)}{\G(2-c_1)\G(c_2)}\\
& &\times e^{-\pi\sqrt{-1}(a+b-c_1-c_2)}x_2^{c_2-1}
F_4(c_2-a+1, c_2-b-1, 2-c_1, c_2;x),\\
f_4^\vee(x)&=&
\frac{x_1^{c_1-1}x_2^{c_2-1}\G(1-c_1)\G(1-c_2)\G(1+b)}{\G(3-c_1-c_2+b)}
F_4(c_1+c_2-a, c_1+c_2-b-2,  c_1, c_2;x).
\end{eqnarray*}
Since $H$ is diagonal, we can easily evaluate $H^{-1}=\tr H^{-1}$.
By multiplying both sides
of (\ref{eq:TPR11}) by $(1-c_1)(1-c_2)/(2\pi\sqrt{-1})^2$ and using the formula $\G(a)\G(1-a)=\pi/\sin(\pi a)$,  
we reduce this relation to the first identity. 
By multiplying the identities arising from the (4,4) and (1,4) entries of (\ref{eq:TPR}) by $(1-c_1)(1-c_2)/(2\pi\sqrt{-1})^2$,
we have the second and third equalities in this corollary.
\qqed

\end{document}